\documentclass[english,a4paper]{article}
\usepackage{mathptmx}
\usepackage[T1]{fontenc}
\usepackage[latin9]{inputenc}
\usepackage{rotfloat}
\usepackage{amsmath}
\usepackage{amsthm}
\usepackage{amssymb}
\usepackage{mathdots}
\usepackage{graphicx}

\makeatletter

\providecommand{\tabularnewline}{\\}

\newcommand{\lyxaddress}[1]{
\par {\raggedright #1
\vspace{1.4em}
\noindent\par}
}
\theoremstyle{plain}
\newtheorem{thm}{\protect\theoremname}
  \theoremstyle{plain}
  \newtheorem{prop}[thm]{\protect\propositionname}
  \theoremstyle{plain}
  \newtheorem{lem}[thm]{\protect\lemmaname}
  \theoremstyle{remark}
  \newtheorem*{rem*}{\protect\remarkname}
  \theoremstyle{plain}
  \newtheorem{cor}[thm]{\protect\corollaryname}

\newcommand{\trace}{\textrm{trace}}
\newcommand{\Real}{\mathbb{R}}
\newcommand{\Comp}{\mathbb{C}}

\makeatother

\usepackage{babel}
  \providecommand{\corollaryname}{Corollary}
  \providecommand{\lemmaname}{Lemma}
  \providecommand{\propositionname}{Proposition}
  \providecommand{\remarkname}{Remark}
\providecommand{\theoremname}{Theorem}

\begin{document}

\title{On a variable step size modification of Hines' method in computational
neuroscience}

\author{Michael Hanke\thanks{email: hanke@nada.kth.se}}
\maketitle

\lyxaddress{KTH Royal Institute of Technology, Department of Mathematics, 100
44 Stockholm, Sweden}
\begin{abstract}
For simulating large networks of neurons Hines proposed a method which
uses extensively the structure of the arising systems of ordinary
differential equations in order to obtain an efficient implementation.
The original method requires constant step sizes and produces the
solution on a staggered grid. In the present paper a one-step modification
of this method is introduced and analyzed with respect to their stability
properties. The new method allows for step size control. Local error
estimators are constructed. The method has been implemented in matlab
and tested using simple Hodgkin-Huxley type models. Comparisons with
standard state-of-the-art solvers are provided.


\emph{Classification}: AMS MSC (2010) 65L20, 65L05, 65L06, 92C42
\end{abstract}

\section{Introduction}

When simulating large networks of neurons a considerable part of the
model consists of the electrical subsystem which in turn extensively
uses the classical Hodkin-Huxley model of nerve activity. Owing to
the large size of the networks to be modeled the efficient numerical
solution of the arising high-dimensional system of ordinary differential
equations is of extraordinary importance. Complex program systems,
e.g., NEURON \cite{CaHi05}, GENESIS \cite{BoBe98}, and many others
are in routine use in order to solve them. In order to construct efficient
numerical methods it is necessary to tailor the methods to the special
properties of the system. The building block is often (variants of)
the Hodgkin-Huxley system \cite{HoHu54}
\begin{align}
C\frac{dV}{dt} & =I(t)-g_{K}n^{4}(V-V_{K})-g_{Na}m^{3}h(V-V_{Na})-g_{L}(V-V_{L}),\label{eq:HHV}\\
\frac{dn}{dt} & =\alpha_{n}(V)(1-n)-\beta_{n}(V)n,\label{eq:HHn}\\
\frac{dm}{dt} & =\alpha_{m}(V)(1-m)-\beta_{m}(V)m,\label{HHm}\\
\frac{dh}{dt} & =\alpha_{h}(V)(1-h)-\beta_{h}(V)h.\label{eq:HHh}
\end{align}
The coefficients $\alpha_{i}(V)$ and $\beta_{i}(V)$ are highly nonlinear
functions of their argument. The key observation in this system is
that the differential equation for the voltage $V$ is linear in $V$
while the differential equations for the gate variables $n,m,h$ are
linear in those. Slightly more general, this system has the structure
\begin{align}
x' & =A(y)x+b(y,t),\label{eq:semilina}\\
y' & =c(x,t)+D(x)y.\label{eq:semilinb}
\end{align}
Here, $x$ denotes the voltage while $y$ is the vector of gate variables,
or vice versa. In general, this model leads to stiff differential
equations such that implicit time stepping methods are necessary.
Standard approaches require the solution of a nonlinear system of
equations in every step by using variants of Newton's method. Given
the large dimension of the usual models, this property may become
a severe restriction. Hines \cite{Hi84} came up with the idea to
discretize the system in two steps: First, the differential equation
for $x$ is discretized leading to a linear system to be solved. Then,
the differential equation for $y$ is discretized. Also here it remains
only a linear system to be solved in contrast to a fully nonlinear
system in the standard approach. Note that the lower dimensional systems
have often a very special structure such that they can be solved very
efficiently. In particular, for the Hodgkin-Huxley system (\ref{eq:HHV})
\textendash{} (\ref{eq:HHh}), both systems have a diagonal system
matrix.

Hines chose the implicit midpoint rule as the basic discretization.
The discrete approximations of $x$ and $y$ are defined on a staggered
grid. In order to fix notation, let $t\in[0,T]$ for some $T>0$ and
$h>0$ be a given step size. For $n=0,1,2,\ldots$ let
\[
t_{n}=nh,\quad t_{n+1/2}=(n+1/2)h.
\]
Hines method reads
\begin{align}
x_{n+1} & =x_{n}+h\left(A(y_{n+1/2})\frac{x_{n+1}+x_{n}}{2}+b(y_{n+1/2},t)\right),\label{eq:hinesa}\\
y_{n+3/2} & =y_{n+1/2}+h\left(c(x_{n+1},t_{n+1})+D(x_{n+1})\frac{y_{n+3/2}+y_{n+1/2}}{2}\right).\label{eq:hinesb}
\end{align}
This method has the following properties:
\begin{itemize}
\item The approximations are available on a staggered grid, only: $x_{n}\approx x(t_{n})$
and $y_{n+1/2}\approx y(t_{n+1/2})$.
\item Since initial values are available for $t=0$ only, the first approximation
$y_{1/2}$ must be computed by other means.
\item The method is second order accurate. However, this property is only
preserved if the step size is constant.
\end{itemize}
In particular the last property calls for a modification of this method
such that a step size control becomes possible.

In this paper, we consider the following modification of Hines' method:\footnote{Gustaf Söderlind (2013), personal communication.}
\begin{align}
x_{n+1/2} & =x_{n}+\frac{h}{2}A(y_{n})x_{n}+b(y_{n},t_{n}),\label{eq:soeda}\\
y_{n+1} & =y_{n}+h\left(c(x_{n+1/2},t_{n+1/2})+D(x_{n+1/2})\frac{y_{n+1}+y_{n}}{2}\right),\label{eq:soedb}\\
x_{n+1} & =x_{n+1/2}+\frac{h}{2}(A(y_{n+1})x_{n+1}+b(y_{n+1},t_{n+1})).\label{eq:soedc}
\end{align}
Since the first of these three equation is an explicit Euler step,
the computational work of the latter method is only slightly more
expensive than that of the original proposal by Hines. However, the
modified version is a genuine one-step method allowing for a step
size change without sacrificing the order.

In this note we will investigate the numerical properties of this
method. In particular, we are interested in the asymptotic stability
of this method. Moreover, we will propose an efficient implementation.
In the final section a few numerical examples will be provided.

\section{Properties of the method}

In order to simplify the notation slightly, we consider the system
\begin{align}
x' & =f(x,y,t),\label{eq:sysa}\\
y' & =g(x,y,t).\label{eq:sysb}
\end{align}
The method (\ref{eq:soeda}) \textendash{} (\ref{eq:soedc}) reduces
to
\begin{align}
x_{n+1/2} & =x_{n}+\frac{h}{2}f(x_{n},y_{n},t_{n}),\label{eq:soedga}\\
y_{n+1} & =y_{n}+hg(x_{n+1/2},\frac{1}{2}(y_{n+1}+y_{n}),t_{n+1/2}),\label{eq:soedgb}\\
x_{n+1} & =x_{n+1/2}+\frac{h}{2}f(x_{n+1},y_{n+1},t_{n+1}).\label{eq:soedgc}
\end{align}
By eliminating the intermediate approximation $x_{n+1/2}$, this system
reduces to
\begin{align}
x_{n+1} & =x_{n}+\frac{h}{2}(f(x_{n+1},y_{n+1},t_{n+1})+f(x_{n},y_{n},t_{n})),\label{eq:soedgsa}\\
y_{n+1} & =y_{n}+hg(x_{n}+\frac{h}{2}f(x_{n},y_{n},t_{n}),\frac{1}{2}(y_{n+1}+y_{n}),t_{n+1/2}).\label{eq:saedgsb}
\end{align}

Let $(x^{\ast},y^{\ast})$ denote the solution of (\ref{eq:sysa})
\textendash{} (\ref{eq:sysb}) on $[0,T]$ subject to the initial
condition $x(0)=x_{0}$ and $y(0)=y_{0}$. Uniqueness is guaranteed
if $f$ and $g$ are Lipschitz continuous with respect to $x$ and
$y$ and continuous with respect to $t$ in a neighborhood $U$ of
the trajectory $\Gamma=\{(x^{\ast}(t),y^{\ast}(t),t)|t\in[0,T]\}$.

Define, for sequences $\{(x_{n},y_{n})\}_{n=0}^{N(h)}$, the discretization
error by
\begin{align*}
{\cal N}_{x}(x_{n},y_{n}) & =\frac{x_{n+1}-x_{n}}{h}-\frac{1}{2}(f(x_{n}y_{n}t_{n})+f(x_{n+1},y_{n+1},t_{n+1})),\\
{\cal N}_{y}(x_{n},y_{n}) & =\frac{y_{n+1}-y_{n}}{h}-g(x_{n}+\frac{h}{2}f(x_{n},y_{n},t_{n}),\frac{1}{2}(y_{n+1}+y_{n}),t_{n+1/2}).
\end{align*}
The method is called stable if there exist an $h_{0}$ and a $K$
such that, for $h<h_{0}$ and for any sequences $\{(x_{n}^{j},y_{n}^{j})\}_{n=0}^{N(h)}$,
$j=1,2$ belonging to $U$ it holds
\begin{multline*}
\max_{n=0,\ldots N(h)}(|x_{n}^{1}-x_{n}^{2}|+|y_{n}^{1}-y_{n}^{2}|)\leq K\left\{ |x_{0}^{1}-x_{0}^{2}|+|y_{0}^{1}-y_{0}^{2}|+\right.\\
\left.\max_{n=0,\ldots,N(h)}\left(|{\cal N}_{x}(x_{n}^{1},y_{n}^{1})-{\cal N}_{x}(x_{n}^{2},y_{n}^{2})|+|{\cal N}_{y}(x_{n}^{1},y_{n}^{1})-{\cal N}_{y}(x_{n}^{2},y_{n}^{2})|\right)\right\} .
\end{multline*}
This definition follows \cite[Section 5.2.3]{AsPe98}.
\begin{prop}
\label{thm:stability}Let $(x^{\ast},y^{\ast})$ be solutions of (\ref{eq:sysa})
\textendash{} (\ref{eq:sysb}) on $[0,T]$ subject to the initial
condition $x(0)=x_{0}$ and $y(0)=y_{0}$. Let $f$ and $g$ be Lipschitz
continuous with respect to $x$ and $y$ and continuous with respect
to $t$ in a neighborhood $U$ of the trajectory $\Gamma=\{(x^{\ast}(t),y^{\ast}(t),t)|t\in[0,T]\}$.
Then, the method (\ref{eq:soedga}) \textendash{} (\ref{eq:soedgc})
is stable.
\end{prop}
\begin{proof}
The method is a combination of two implicit Runge-Kutta methods. So
the proof is standard.
\end{proof}
\begin{prop}
Let the assumptions of Theorem~\ref{thm:stability} be fulfilled
and $f$ and $g$ be sufficiently often differentiable. Then, (\ref{eq:soedga})
\textendash{} (\ref{eq:soedgc}) is a second order method and the
global error has an asymptotic expansion in powers of $h^{2}$.
\end{prop}
\begin{proof}
As a Runge-Kutta method, the local error possesses an expansion in
powers of $h$, \cite[Theorem 3.2]{HaNoWa92}. Taylor expansion shows
that, for the truncation error, it holds
\[
{\cal N}_{x}(x^{\ast}(t_{n}),y^{\ast}(t_{n}))=O(h^{2}),\quad{\cal N}_{y}(x^{\ast}(t_{n}),y^{\ast}(t_{n}))=O(h^{2}).
\]
 Since the method is stable, it is second order convergent. 

Rewriting (\ref{eq:soedgsa}), it holds $x_{n}+(h/2)f(x_{n},y_{n},t_{n})=x_{n+1}-(h/2)f(x_{n+1},y_{n+1},t_{n+1})$
such that (\ref{eq:saedgsb}) can be rewritten as
\[
y_{n}=y_{n+1}-hg(x_{n+1}-(h/2)f(x_{n+1},y_{n+1},t_{n+1}),\frac{1}{2}(y_{n+1}+y_{n}),t_{n+1/2}).
\]
Hence, the method is symmetric and the assertion about the asymptotic
expansion follows from \cite[Theorem 8.10]{HaNoWa92}.
\end{proof}
A more interesting question is the asymptotic properties of this method.
Recall that the system to be solved is usually stiff. This is also
the case if the two components are considered individually, that is
(\ref{eq:sysa}) for fixed $y$ and (\ref{eq:sysb}) for fixed $x$.
So the standard notions of asymptotic stability that base on the test
scalar equation $z'=\lambda z$ are not useful here. A more appropriate
test equation must have at least two components. This leads to the
proposal
\begin{align}
x' & =\mu x+ay,\label{eq:testa}\\
y' & =bx+\lambda y,\label{eq:testb}
\end{align}
 where
\begin{equation}
\mu,\lambda<0\text{ and }ab<\mu\lambda.\label{eq:Stabi}
\end{equation}
Under these conditions, the system is asymptotically stable as well
as the individual components are. This test system has been proposed
by Strehmel\&Weiner \cite{StWe84} in order to characterize stability
properties of partitioned Runge-Kutta methods.

The method (\ref{eq:soedga}) \textendash{} (\ref{eq:soedgc}) applied
to (\ref{eq:testa}) \textendash{} (\ref{eq:testb}) gives rise to
\begin{equation}
A\left(\begin{array}{c}
x_{n+1}\\
y_{n+1}
\end{array}\right)=B\left(\begin{array}{c}
x_{n}\\
y_{n}
\end{array}\right)\label{eq:s}
\end{equation}
where
\begin{align*}
A & =\left(\begin{array}{cc}
1-\frac{h\mu}{2} & -\frac{ha}{2}\\
0 & 1-\frac{h\lambda}{2}
\end{array}\right),\\
B & =\left(\begin{array}{cc}
1+\frac{h\mu}{2} & \frac{ha}{2}\\
hb\left(1+\frac{h\mu}{2}\right) & 1+\frac{h\lambda}{2}+\frac{abh^{2}}{2}
\end{array}\right).
\end{align*}
This recursion can be rewritten in the form
\[
\left(\begin{array}{c}
x_{n+1}\\
y_{n+1}
\end{array}\right)=C\left(\begin{array}{c}
x_{n}\\
y_{n}
\end{array}\right)
\]
with
\[
C=\left(\begin{array}{cc}
\alpha\left(1+\gamma\frac{h\mu}{2}(\beta-1)\right) & ha\left(\frac{1}{\left(1-\frac{h\mu}{2}\right)\left(1-\frac{h\lambda}{2}\right)}+\frac{\gamma}{4}(\alpha-1)(\beta-1)\right)\\
hb\frac{1+\frac{h\mu}{2}}{1-\frac{h\lambda}{2}} & \beta+\gamma(\beta-1)\frac{h\mu}{2}
\end{array}\right)
\]
where
\[
\alpha=\frac{1+\frac{h\mu}{2}}{1-\frac{h\mu}{2}},\quad\beta=\frac{1+\frac{h\lambda}{2}}{1-\frac{h\lambda}{2}}
\]
are the stability functions of the midpoint and trapezoidal rule,
respectively, applied to the equations (\ref{eq:testa}) \textendash{}
(\ref{eq:testb}) individually. Note that, under the conditions (\ref{eq:Stabi}),
it holds $|\alpha|<1$ and $|\beta|<1$. 

The recursion is asymptotically stable if and only if the eigenvalues
of $C$ are less than one in absolute value. A short computation shows
that the characteristic polynomial of $C$ has the form
\[
\chi(s)=s^{2}-(\alpha+\beta+\gamma(\alpha-1)(\beta-1))s+\alpha\beta
\]
were $\gamma=\frac{ab}{\mu\lambda}$.
\begin{lem}
\label{lem:quadrat}Consider the polynomial $\chi(s)=s^{2}-(\alpha+\beta+\gamma(\alpha-1)(\beta-1))s+\alpha\beta$
with $|\alpha|<1$ and $|\beta|<1$. For the roots $s_{1}$ and $s_{2}$
of $\chi(s)=0$ it holds $\max(|s_{1}|,|s_{2}|)<1$ if and only if
\[
-\frac{(1+\alpha)(1+\beta)}{(1-\alpha)(1-\beta)}<\gamma<1.
\]
\end{lem}
\begin{rem*}
Under the assumption (\ref{eq:Stabi}) it holds $\gamma<1$. Note
that $\gamma$ may be negative.
\end{rem*}
\begin{proof}
We use the change of variables
\[
w(z)=\frac{1+z}{1-z}.
\]
$w$ maps the negative complex halfplane $\Comp^{-}=\{z\in\Comp|\mathfrak{R}z<0\}$
uniquely onto the disk $S=\{z\in\Comp||z|<1\}$. So it holds $s_{1},s_{2}\in S$
if and only if $z_{1},z_{2}\in\Comp^{-}$ for the zeros of $\chi(w(z))$.
Denote for short $\chi(s)=s^{2}+a_{1}s+a_{0}$. A short calculation
provides
\[
\chi(w(z))=\frac{z^{2}(1-a_{1}+a_{0})+z(2-2a_{0})+(1+a_{1}+a_{0})}{(1-z)^{2}}.
\]
The zeros of this functions are those of the enumerator polynomial
$\eta(z)=c_{2}z^{2}+c_{1}z+c_{0}$ where
\begin{align*}
c_{2} & =1+(\alpha+\beta+\gamma(\alpha-1)(\beta-1))+\alpha\beta,\\
c_{1} & =2-2\alpha\beta,\\
c_{0} & =1-(\alpha+\beta+\gamma(\alpha-1)(\beta-1))+\alpha\beta.
\end{align*}
According to the Routh-Hurwitz criterion \cite[Theorem I.13.4]{HaNoWa92}
the roots of $\eta$ lie all in $\Comp^{-}$ if and only if all coefficients
$c_{i}$ have the same sign. Since $|\alpha\beta|<1$, it holds $c_{1}>0$.
The condition $c_{0}>0$ is equivalent to 
\[
\frac{1+\alpha\beta-\alpha-\beta}{(\alpha-1)(\beta-1)}=1>\gamma
\]
 while $c_{2}>0$ is equivalent to
\[
\frac{-(\alpha+1)(\beta+1)}{(\alpha-1)(\beta-1)}<\gamma.
\]
This completes the proof.
\end{proof}
\begin{thm}
Under the assumption $\mu,\lambda<0$, the recursion (\ref{eq:s})
is asymptotically stable if and only if
\[
-\frac{(1+\alpha)(1+\beta)}{(1-\alpha)(1-\beta)}<\gamma<1.
\]
\end{thm}
The assertion is a consequence of Lemma~\ref{lem:quadrat}.
\begin{cor}
For every $\gamma<0$, there exists a $h(\mu,\lambda,\gamma)>0$ such
that the recursion (\ref{eq:s}) is unstable for all $h>h(\mu,\lambda,\gamma)$.
\end{cor}
\begin{proof}
One can easily show that $\alpha=\alpha(h)$ and $\beta=\beta(h)$
are monotonically decreasing functions of $h$. Moreover,
\[
\varphi(h)=\frac{1+\alpha}{1-\alpha}
\]
 is a monotonically decreasing function of $\alpha$ and it holds
$\lim_{h\rightarrow0}\varphi(\alpha(h))=\infty$ and $\lim_{h\rightarrow\infty}\varphi(\alpha(h))=0$.
\end{proof}
\begin{rem*}
It is interesting to see how the corresponding results for the original
Hines' method look like. The recursion becomes
\[
\left(\begin{array}{c}
x_{n+1}\\
y_{n+3/2}
\end{array}\right)=C_{\textrm{Hines}}\left(\begin{array}{c}
x_{n}\\
y_{n+1/2}
\end{array}\right)
\]
with
\begin{align*}
C_{\textrm{Hines}} & =\left(\begin{array}{cc}
\alpha & \frac{ha}{\left(1-\frac{h\mu}{2}\right)}\\
\alpha\frac{hb}{\left(1-\frac{h\lambda}{2}\right)} & \beta+\frac{abh^{2}}{\left(1-\frac{h\mu}{2}\right)\left(1-\frac{h\lambda}{2}\right)}
\end{array}\right)\\
 & =\left(\begin{array}{cc}
\alpha & \frac{ha}{\left(1-\frac{h\mu}{2}\right)}\\
\alpha\frac{hb}{\left(1-\frac{h\lambda}{2}\right)} & \beta+\gamma(1-\alpha)(1-\beta)
\end{array}\right).
\end{align*}
It holds
\begin{align*}
\trace(C_{\textrm{Hines}}) & =\alpha+\beta+\gamma(1-\alpha)(1-\beta),\\
\det(C_{\textrm{Hines}}) & =\alpha\beta+\alpha\gamma(1-\alpha)(1-\beta)-\alpha\frac{abh^{2}}{\left(1-\frac{h\mu}{2}\right)\left(1-\frac{h\lambda}{2}\right)}\\
 & =\alpha\beta.
\end{align*}
So the stability polynomial of Hines' method and its modification
(\ref{eq:soeda}) \textendash{} (\ref{eq:soedc}) are identical.
\end{rem*}

\section{The relation of this method to Strang's splitting and the Peaceman-Rachford
method}

\subsection{Strang's splitting}

In this section, we will consider an approach to solving (\ref{eq:sysa})
\textendash{} (\ref{eq:sysb}) using a splitting method. For that,
let
\[
U=\left(\begin{array}{c}
x\\
y
\end{array}\right),\quad F(U,t)=\left(\begin{array}{c}
f(x,y,t)\\
g(x,y,t)
\end{array}\right).
\]
Then (\ref{eq:sysa}) \textendash{} (\ref{eq:sysb}) is equivalent
to $U'=F(U,t)$. Introduce the splitting
\begin{equation}
F(U,t)=F_{1}(U,t)+F_{2}(U,t)=\left(\begin{array}{c}
f(x,y,t)\\
0
\end{array}\right)+\left(\begin{array}{c}
0\\
g(x,y,t)
\end{array}\right).\label{eq:decomp}
\end{equation}
A classical example of operator splitting is Strang's approach \cite{St68}.
In order to advance the solution one step of step size $h$ from $t_{n}$
to $t_{n+1}$ the system $U'=F_{1}(U,t)$ is first integrated over
the interval $[t_{n},t_{n}+h/2]$, then the results are used to integrate
the system $U'=F_{2}(U,t)$ on $[t_{n},t_{n+1}]$, and finally the
first system on $[t_{n+1}-h/2,t_{n+1}]$. For the decomposition (\ref{eq:decomp}),
this gives rise to the following three steps:
\begin{enumerate}
\item Integrate $x'=f(x,y_{n},t)$, $x(t_{n})=x_{n}$ on $[t_{n,}t_{n}+h/2]$.
Denote $x_{n+1/2}=x(t_{n}+h/2)$.
\item Integrate $y'=g(x_{n+1/2},y,t)$, $y(t_{n})=y_{n}$ on $[t_{n},t_{n+1}]$.
Denote $y_{n+1}=y(t_{n+1})$.
\item Integrate $x'=f(x,y_{n+1},t)$, $x(t_{n}+h/2)=x_{n+1/2}$ on $[t_{n+1}-h/2,t_{n+1}]$.
Let $x_{n+1}=x(t_{n+1})$.
\end{enumerate}
Even if $f,g$ are linear functions, the operators $F_{1}$ and $F_{2}$
do not commute. So we expect the splitting to be second order accurate
at least in the linear case (e.g., \cite[Chapter 4]{HuVe03}). A comparison
with the finite difference method (\ref{eq:soedga}) \textendash{}
(\ref{eq:soedgc}) reveals that the latter can be interpreted as a
second order discretization of Strang's splitting. 

Let us ask the question of asymptotic stability of Strang's splitting
applied to the test system (\ref{eq:testa}) \textendash{} (\ref{eq:testb}).
Similarly as before the recursion can be written down in the form
\begin{align}
\left(\begin{array}{c}
x_{n+1}\\
y_{n+1}
\end{array}\right) & =C_{\text{Strang}}\left(\begin{array}{c}
x_{n}\\
y_{n}
\end{array}\right)\label{eq:ss}\\
C_{\text{Strang}} & =\left(\begin{array}{cc}
\alpha+\alpha^{1/2}T & \frac{a}{\mu}\left(\alpha^{1/2}-1\right)\left(\alpha^{1/2}+T+\beta\right)\\
\frac{b}{\lambda}\alpha^{1/2}\left(\beta-1\right) & T+\beta
\end{array}\right)\nonumber 
\end{align}
where $T=\gamma\left(\alpha^{1/2}-1\right)\left(\beta-1\right)$.
Here,
\[
\alpha=e^{\mu h},\quad\beta=e^{\lambda h}.
\]
The characteristic polynomial becomes
\[
\chi(s)=s^{2}-(\alpha+\beta+\gamma(\alpha-1)(\beta-1))s+\alpha\beta.
\]
This is structurally identical to the one for the discrete case.
\begin{thm}
Let $\mu<0$ and $\lambda<0$. The recursion (\ref{eq:ss}) is asymptotically
stable if and only if 
\[
-\frac{\left(1+e^{h\mu}\right)\left(1+e^{h\lambda}\right)}{\left(e^{h\mu}-1\right)\left(e^{h\lambda}-1\right)}<\gamma<1.
\]
\end{thm}
\begin{proof}
The result is a consequence of the results of Lemma~\ref{lem:quadrat}
by setting $\alpha=e^{h\mu}$ and $\beta=e^{h\lambda}$.
\end{proof}
Note that always $\gamma<1$ under the assumption $ab<\mu\lambda$.
Moreover, for any $h>0,\mu<0,\lambda<0$ and
\[
\psi(h)=\frac{\left(1+e^{h\mu}\right)\left(1+e^{h\lambda}\right)}{\left(e^{h\mu}-1\right)\left(e^{h\lambda}-1\right)}
\]
it holds that $\psi$ is a monotonically decreasing function with
$\lim_{h\rightarrow0}\psi(h)=\infty$ and $\lim_{h\rightarrow\infty}\psi(h)=1$.
We have immediately
\begin{cor}
If $\gamma<-1$ and $\mu<0,\lambda<0$, the there exist always an
$h(\gamma,\mu,\lambda)>0$ such that the recursion (\ref{eq:ss})
is unstable for $h>h(\gamma,\mu,\lambda)$.
\end{cor}
Compared to the discrete case, the stability domain is slightly larger
for Strang's splitting. However, even here the stability domain is
bounded.

\subsection{The Peaceman-Rachford method}

The Peaceman-Rachford method was introduced in \cite{PeRa55} in order
to solve semidiscretized linear parabolic partial differential equations.
If the system (\ref{eq:sysa}) \textendash{} (\ref{eq:sysb}) is splitted
according to (\ref{eq:decomp}), the solution $U$ is then advanced
from $t_{n}$ to $t_{n+1}$ by
\begin{align*}
U_{n+1/2} & =U_{n}+\frac{h}{2}\left(F_{1}(U_{n+1/2},t_{n+1/2})+F_{2}(U_{n},t_{n})\right),\\
U_{n+1} & =U_{n+1/2}+\frac{h}{2}\left(F_{1}(U_{n+1/2},t_{n+1/2})+F_{2}(U_{n+1},t_{n+1})\right).
\end{align*}
The application of the Peaceman-Rachford method in (\ref{eq:decomp})
leads to
\begin{align*}
\left(\begin{array}{c}
x_{n+1/2}\\
y_{n+1/2}
\end{array}\right) & =\left(\begin{array}{c}
x_{n}\\
y_{n}
\end{array}\right)+\frac{h}{2}\left(\left(\begin{array}{c}
0\\
g(x_{n+1/2},y_{n+1/2},t_{n+1/2}
\end{array}\right)+\left(\begin{array}{c}
f(x_{n},y_{n},t_{n})\\
0
\end{array}\right)\right),\\
\left(\begin{array}{c}
x_{n+1}\\
y_{n+1}
\end{array}\right) & =\left(\begin{array}{c}
x_{n+1/2}\\
y_{n+1/2}
\end{array}\right)+\frac{h}{2}\left(\left(\begin{array}{c}
0\\
g(x_{n+1/2},y_{n+1/2},t_{n+1/2}
\end{array}\right)+\left(\begin{array}{c}
f(x_{n+1},y_{n+1},t_{n+1})\\
0
\end{array}\right)\right).
\end{align*}
 Writing out the components, this recursion becomes
\begin{description}
\item [{(i)}] $x_{n+1/2}=x_{n}+\frac{h}{2}f(x_{n},y_{n},t_{n})$
\item [{(ii)}] $y_{n+1/2}=y_{n}+\frac{h}{2}g(x_{n+1/2},y_{n+1/2},t_{n+1/2})$
\item [{(iii)}] $y_{n+1}=y_{n+1/2}+\frac{h}{2}g(x_{n+1/2},y_{n+1/2},t_{n+1/2})$
\item [{(iv)}] $x_{n+1}=x_{n+1/2}+\frac{h}{2}f(x_{n+1},y_{n+1},t_{n+1})$
\end{description}
Inserting (i) into (iv) we obtain
\[
x_{n+1}=x_{n}+\frac{h}{2}\left(f(x_{n},y_{n},t_{n})+f(x_{n+1},y_{n+1},t_{n+1})\right).
\]
Similarly, from (ii) and (iii) we have
\[
y_{n+1}=y_{n}+hg(x_{n}+\frac{h}{2}f(x_{n},y_{n},t_{n}),y_{n+1/2},t_{n+1/2}).
\]
By subtracting (ii) and (iii) we arrive at
\[
y_{n+1/2}=\frac{1}{2}\left(y_{n}+y_{n+1}\right).
\]
Hence, the Peaceman-Rachford method applied to our system becomes
\begin{align*}
x_{n+1} & =x_{n}+\frac{h}{2}\left(f(x_{n},y_{n},t_{n})+f(x_{n+1},y_{n+1},t_{n+1})\right),\\
y_{n+1} & =y_{n}+hg(x_{n}+\frac{h}{2}f(x_{n},y_{n},t_{n}),\frac{1}{2}\left(y_{n}+y_{n+1}\right),t_{n+1/2}).
\end{align*}
So this method is equivalent to (\ref{eq:soedgsa}) \textendash{}
(\ref{eq:saedgsb}).

The Peaceman-Rachford method has been used extensively to solve partial
differential equations. There, one is mainly interested in showing
stability properties independent of the spatial discretization. These
stability estimates are often based on monotonicity assumptions (one-sided
Lipschitz conditions) on the right-hand side $F_{1},F_{2}$. In particular,
let the condition
\[
\left\langle F_{i}(\tilde{w},t)-F_{i}(w,t),\tilde{w}-w\right\rangle \leq\nu\|\tilde{w}-w\|^{2},\quad i=1,2,
\]
 hold for all $\tilde{w},w$ and $t$ with a certain constant $\nu\in\Real$.
Here, $\left\langle \cdot\right\rangle $ denotes the Euclidean inner
product and $\|\cdot\|$ the Euclidean norm. Hundsdorfer\&Verwer \cite{HuVe89}
show that the method is unconditionally stable (that is, for all step
sizes $h$) if $\nu\leq0$. However, if $\nu>0$, stability can only
be guaranteed if the step size is restricted by $h\nu<2$.

How do these results translate to our model system
\begin{align*}
x' & =\mu x+ay,\\
y' & =bx+\lambda y,
\end{align*}
where $\mu,\lambda<0$ and $ab<\mu\lambda$?

In the linear autonomous case, the monotonicity requirement reduces
to
\[
w^{T}F_{i}(w)\leq\nu\|w\|^{2}\quad\text{for all }w.
\]
We have
\begin{align*}
w^{T}F_{2}(w)\leq\nu\|w\|^{2}\text{ for all }w & \Longleftrightarrow y(bx+\lambda y)\leq\nu(x^{2}+y^{2})\text{ for all }x,y\\
 & \Longleftrightarrow0\leq\nu x^{2}-bxy+(\nu-\lambda)y^{2}\text{ for all }x,y.
\end{align*}
For the latter inequality to hold for all $x,y$, we must have $\nu\geq0$.
\begin{align*}
w^{T}F_{2}(w)\leq\nu\|w\|^{2}\text{ for all }w & \Longleftrightarrow0\leq\left(\sqrt{\nu}x-\frac{1}{2\sqrt{\nu}}by\right)^{2}+\left(\nu-\lambda-\frac{b^{2}}{4\nu}\right)y^{2}\text{ for all }x,y\\
 & \Longleftrightarrow0\leq\nu-\lambda-\frac{b^{2}}{4\nu}\\
 & \Longleftrightarrow\frac{1}{2}(\lambda+\sqrt{\lambda^{2}+b^{2}})\leq\nu.
\end{align*}
Similarly, 
\[
w^{T}F_{1}(w)\leq\nu\|w\|^{2}\text{ for all }w\Longleftrightarrow\frac{1}{2}(\mu+\sqrt{\mu^{2}+a^{2}})\leq\nu.
\]
Hence, we have always $\nu>0$ unless $a=b=0$. 

\section{Implementation}

For an efficient implementation, estimations of the local error and
step size control are of utmost importance. In this section we will
discuss these issues.

The method (\ref{eq:soedga}) \textendash{} (\ref{eq:soedgc}) does
not have an imbedded error estimator. So we have two possibilities
for estimating the local error:
\begin{enumerate}
\item Since the discrete solution possesses an asymptotic expansion in powers
of $h^{2}$, the error can be estimated via Richardson extrapolation.
This can be combined with local extrapolation.
\item Use the detailed representation of the leading error term.
\end{enumerate}
The first idea seems to be rather straightforward. However, one has
to keep in mind that the individual equations in (\ref{eq:semilina})
\textendash{} (\ref{eq:semilinb}) are often stiff, and the discretization
reduces to the trapezoidal rule and the implicit midpoint rule, respectively,
for decoupled systems. For such systems and these discretizations,
we expect a simple step size halving to behave very badly when using
local extrapolation since the resulting method has a bounded stability
domain \cite[p. 133]{HaWa96}. In fact, when applying our method and
step size halving to Hodkin-Huxley systems, we observed a behavior
of the method which is typical for instabilities due to too large
step sizes at low tolerances. Therefore, we implemented two versions:
\begin{itemize}
\item a version using the step size subdivisions $\{1,2\}$ without extrapolation
(\texttt{modhines});
\item a version using the step size subdivisions $\{1,3\}$ with local extrapolation
(\texttt{modhext}). 
\end{itemize}
It should be mentioned that a theoretical analysis of the extrapolation
method is missing so far. It is known that the domains of absolute
stability for the extrapolated trapezoidal rule become smaller and
smaller with the number of extrapolation steps. However, the present
method should be investigated using the test system (\ref{eq:testa})
\textendash{} (\ref{eq:testb}).

In practice, there is no problem to further extrapolate in the extrapolation
tableau. However, in the applications we are aiming at, we expect
mainly low accuracies to be required such that high order methods
will not provide much benefit.

The discrete solution $x_{n+1}$ is computed using the trapezoidal
rule (\ref{eq:soedgsa}). Therefore, the local discretization error
has the representation
\[
\tau_{x,n}={\cal N}_{x}(x(t_{n}),y(t_{n}))=-\frac{1}{12}x'''(t_{n})h^{2}+O(h^{3}).
\]
A similar computation leads to
\begin{align*}
\tau_{y,n} & ={\cal N}_{y}(x(t_{n}),y(t_{n}))\\
 & =\left(\frac{1}{24}y'''(t_{n})+\frac{1}{8}\left(\frac{\partial}{\partial x}g(x(t_{n}),y(t_{n}),t_{n})x''(t_{n})-\frac{\partial}{\partial y}g(x(t_{n}),y(t_{n}),t_{n})y''(t_{n})\right)\right)h^{2}+O(h^{3}).
\end{align*}

In order to approximate the derivatives $x'''$, $y'''$, and $y''$,
the discrete solution is interpolated by a 3rd order Hermite polynomial
using the the values $(x_{n},y_{n}),(x_{n+1},y_{n+1})$ and the function
values $(f(x_{n},y_{n},t_{n}),g(x_{n},y_{n},t_{n}))$, $(f(x_{n+1},y_{n+1},t_{n+1}),g(x_{n+1},y_{n+1},t_{n+1}))$.
In case of an accepted step, $f(x_{n},y_{n},t_{n})$ and $f(x_{n+1},y_{n+1},t_{n+1})$
are available for free. It should be noted that in the case of a constant
coefficient system $y'=Ay$ it holds $y'''=Ay''=(\partial/\partial y)gy''$
such that $\tau_{y,n}=-\frac{1}{12}y'''(t_{n})h^{2}+O(h^{3})$. In
our implementation (\texttt{modhnew}) we use the ``simplified''
approximation $\tau_{y,n}=-\frac{1}{12}y'''(t_{n})h^{2}+O(h^{3})$
for $y$ even in the general case.

We tested also a number of step size selection strategies following
proposals in \cite{So02,So03}. The PI.4.2 controller \cite{So02}
was finally chosen.

\section{Numerical examples}

\begin{sidewaysfigure}
\begin{centering}
\begin{tabular}{cc}
\includegraphics[width=0.4\textwidth]{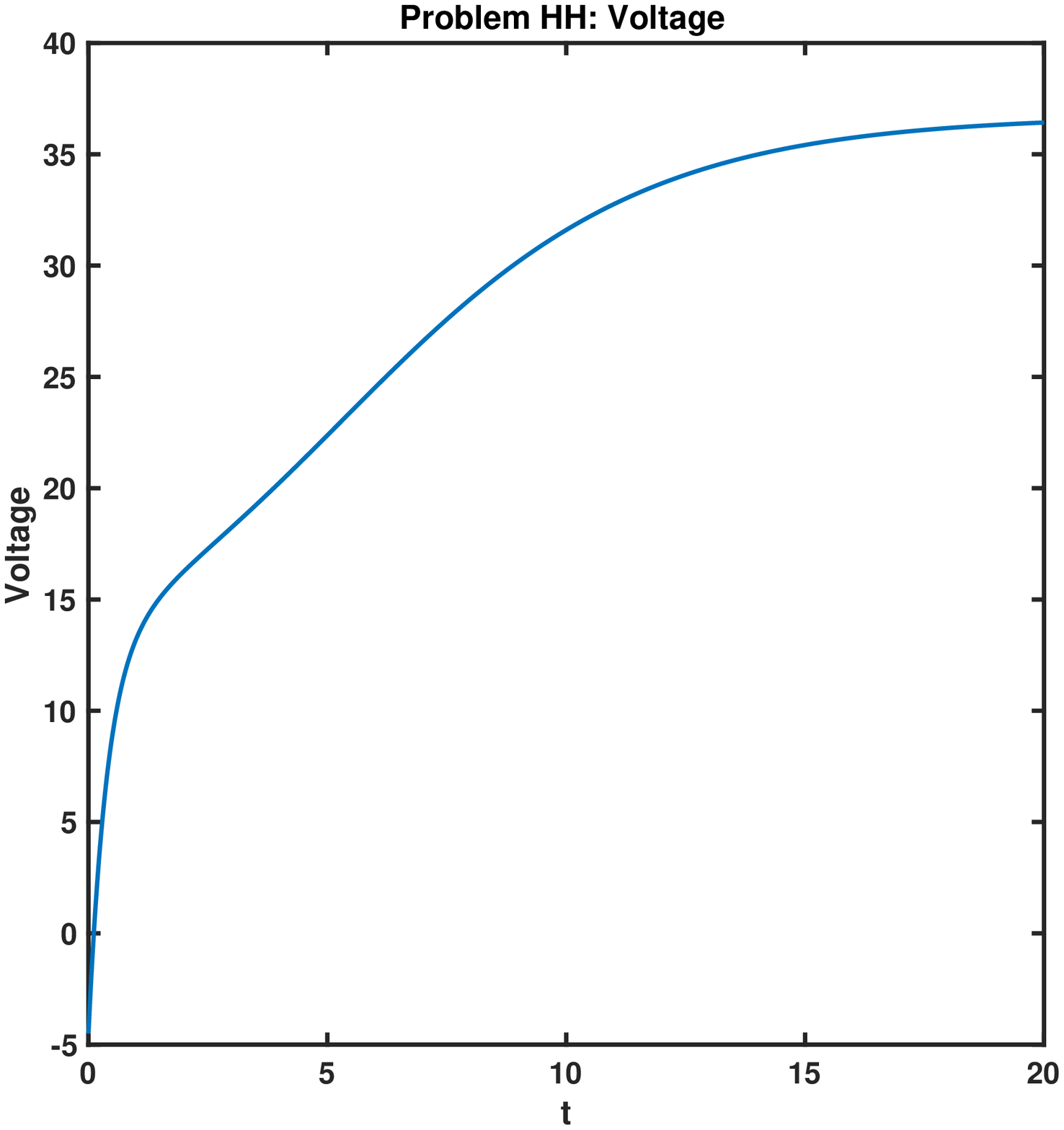} & \includegraphics[width=0.39\textwidth]{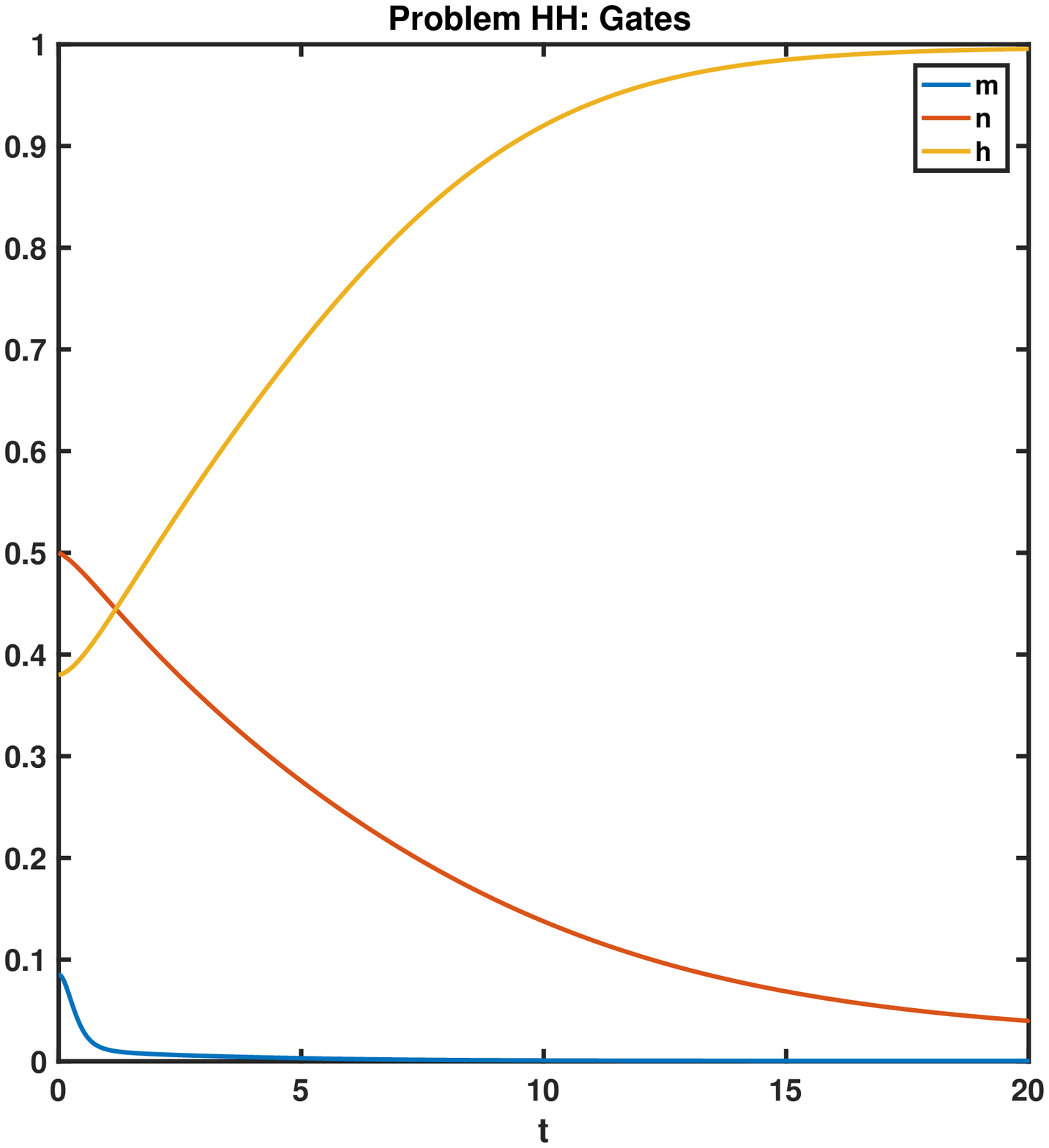}\tabularnewline
(a) & (b)\tabularnewline
\end{tabular}
\par\end{centering}
\caption{Solution of the Hodkin-Huxley system. (a) Voltage $V$, (b) gate variables
$m,n,h$\label{fig:Solution-HH}}
\end{sidewaysfigure}

\begin{sidewaysfigure}
\begin{centering}
\begin{tabular}{ccc}
\includegraphics[width=0.3\textwidth]{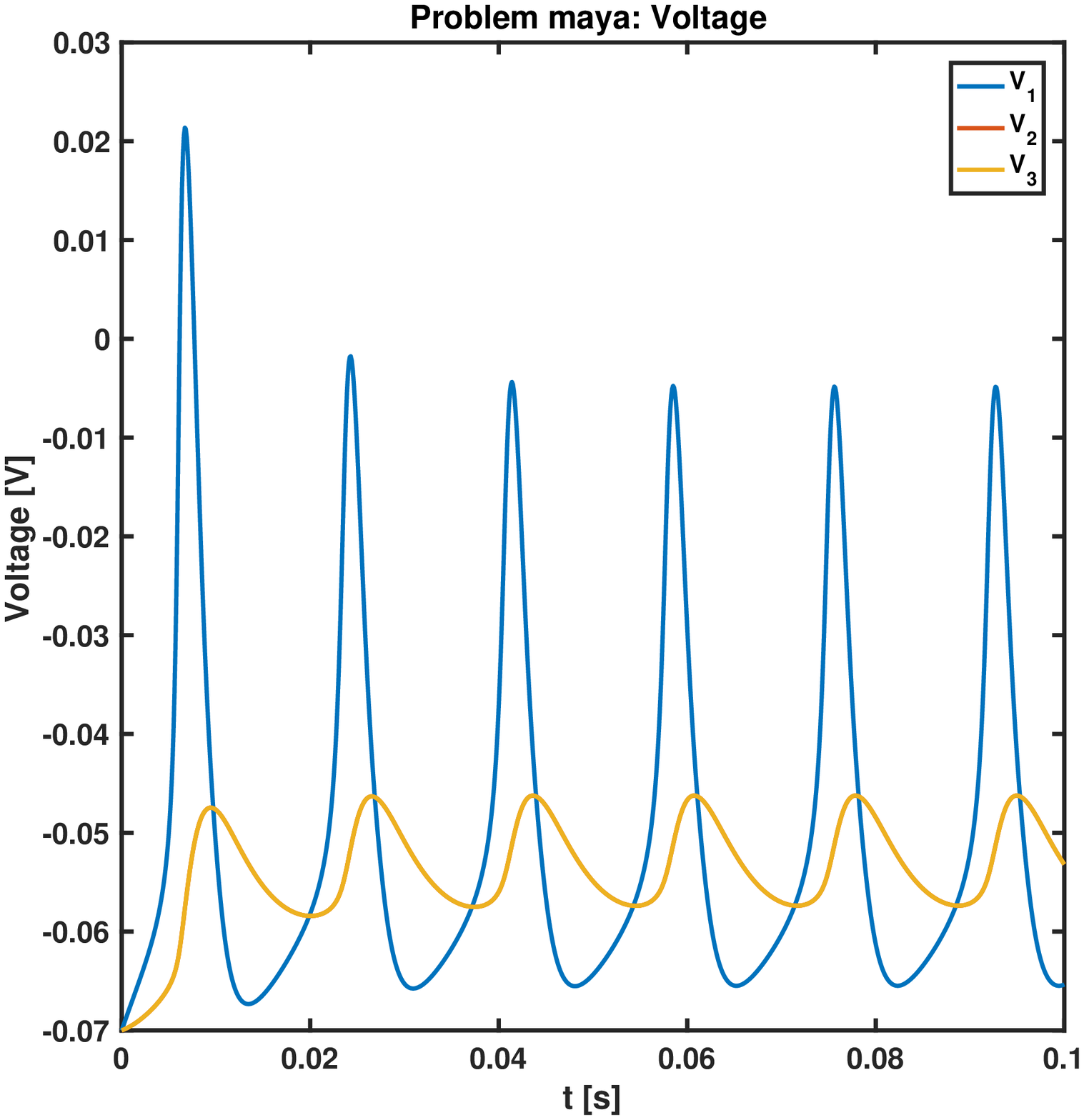} & \includegraphics[width=0.3\textwidth]{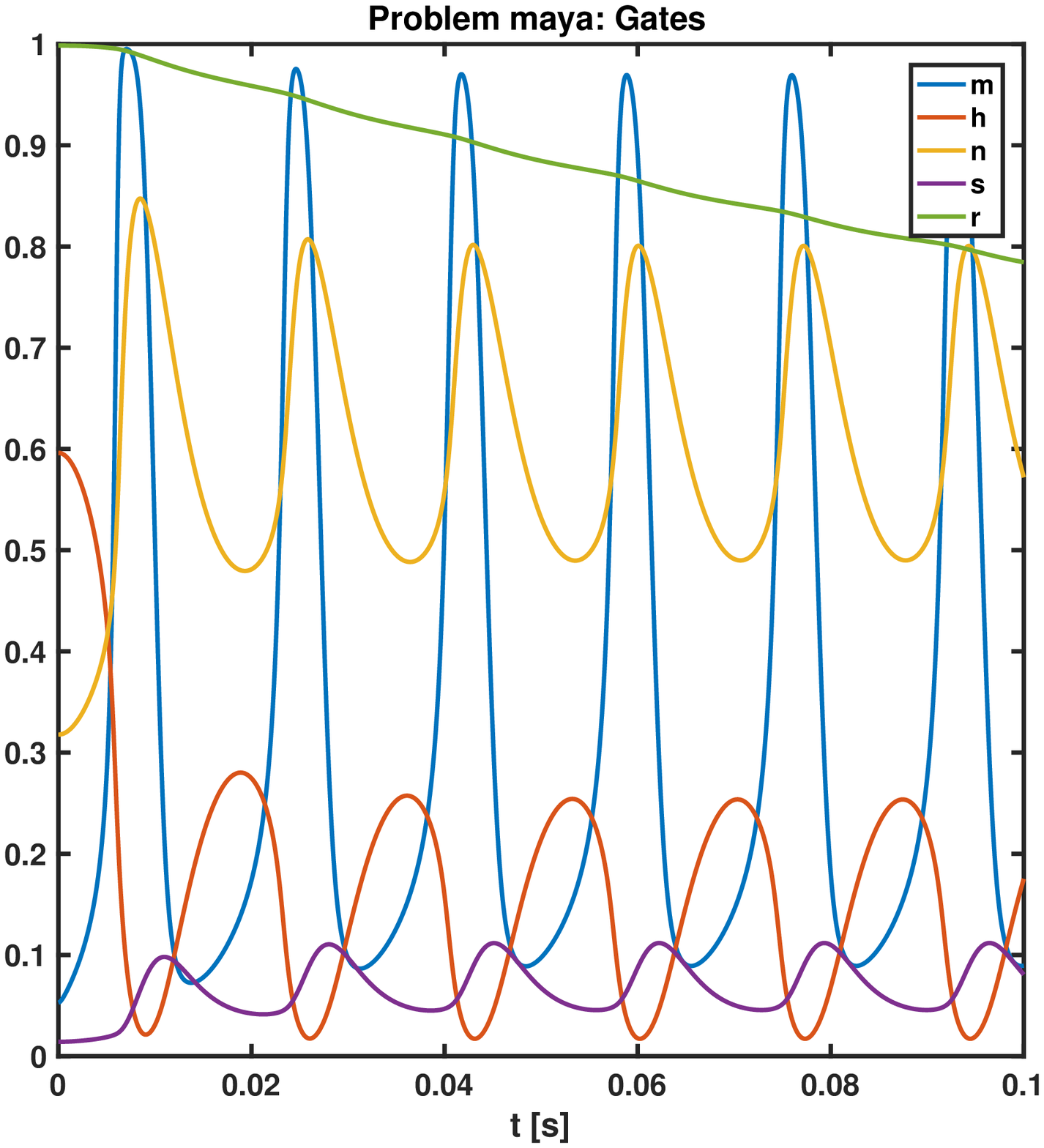} & \includegraphics[width=0.3\textwidth]{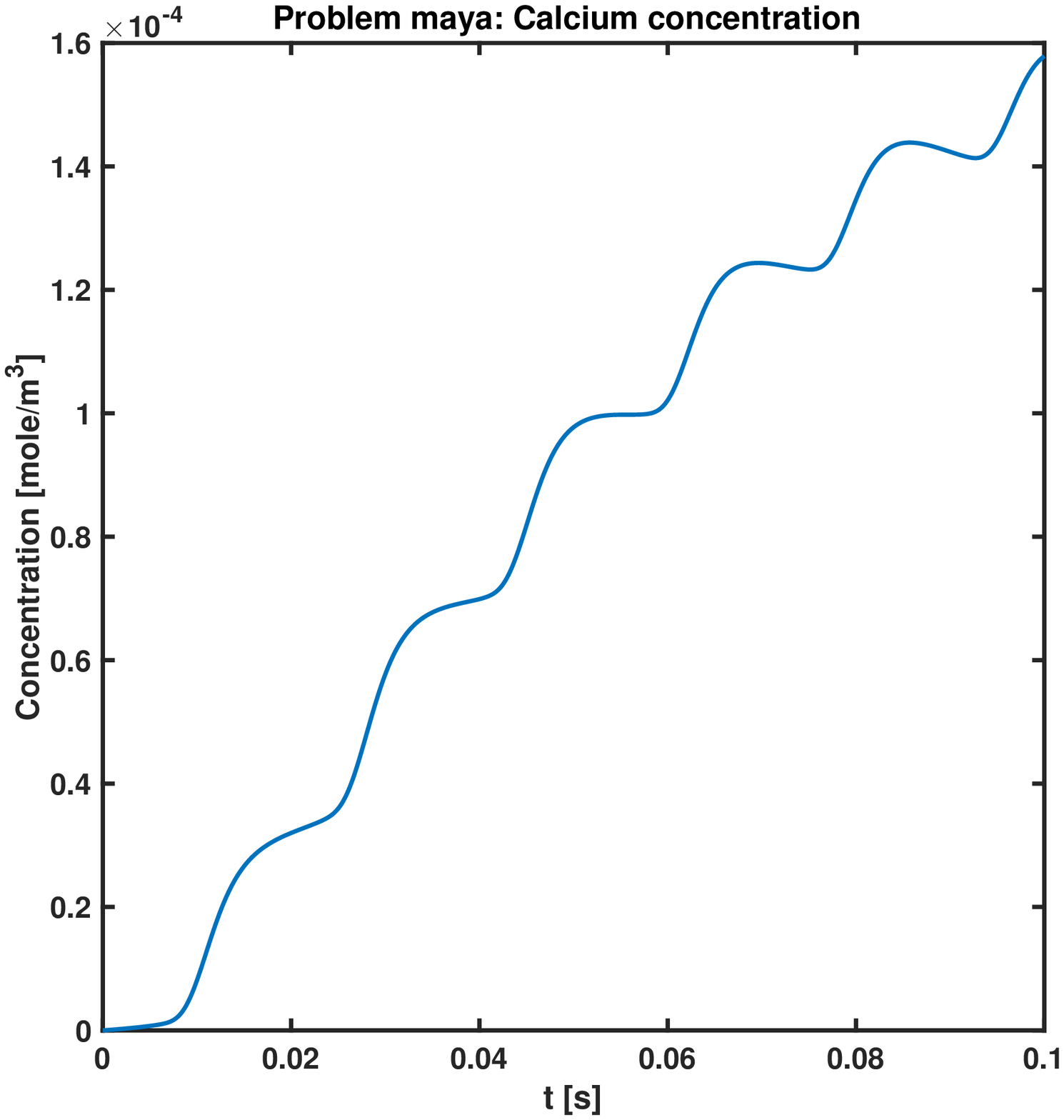}\tabularnewline
(a) & (b) & (c)\tabularnewline
\end{tabular}
\par\end{centering}
\caption{Solution of the soma-dendrite-spine system. (a) Voltages $V_{1},V_{2},V_{3}$.
Note that $V_{2}=V_{3}$ up to plotting accuracy, (b) gate variables
$m,n,h,r,s$, (c) calcium concentration $c_{Ca}$ \label{fig:Solution-of-maya}}
\end{sidewaysfigure}

We show the performance on two benchmark problems; a pure Hodgkin-Huxley
system (Example 1) and a model of a nerve cell with a spine by using
compartment modeling (Example 2).
\begin{description}
\item [{Example~1}] In the Hodkin-Huxley system we choose the parameters
from \cite{HoHu54}:
\[
\begin{gathered}C=1,\quad I(t)=14.2\\
g_{K}=36,\quad g_{Na}=120,\quad g_{L}=0.3\\
V_{K}=12,\quad V_{Na}=-115,\quad V_{L}=-10.599\\
\alpha_{n}(V)=0.1\psi(0.1(V+10)),\quad\beta_{n}(V)=0.125\exp(V/80)\\
\alpha_{m}(V)=\psi(0.1(V+25)),\quad\beta_{m}(V)=4\exp(V/18)\\
\alpha_{h}(V)=0.07\exp(0.05V),\quad\beta_{h}(V)=\left(1+\exp(0.1(V+30))\right)^{-1}\\
\psi(x)=\frac{x}{e^{x}-1}
\end{gathered}
\]
This system is solved on $[0,20]$ with the initial values
\[
V(0)=-4.5,\quad m(0)=0.085,\quad n(0)=0.5,\quad h(0)=0.38.
\]
Figure~\ref{fig:Solution-HH} shows a plot of the solutions.
\item [{Example~2}] We consider the neuron model consisting of three compartments,
the soma, the dendrite and a spine that has been proposed in \cite{Br11}.
Each compartment carries its own potential $V_{i}$, $i=1,2,3$ corresponding
to the soma, the dendrite and the spine. The soma is modeled including
three types of channels ($n,m,h$) while the spine has two of them
($r,s$). There are no channels attached to the dendrite. Additionally,
the calcium ion dynamics is taken into account. The latter contains
an additional degradation in a calcium pool. The equations become
\begin{align*}
C_{1}\frac{dV_{1}}{dt} & =I(t)-g_{K}n^{4}(V_{1}-V_{K})-g_{Na}m^{3}h(V_{1}-V_{Na})+\frac{V_{2}-V_{1}}{R_{a,2}}-\frac{V_{1}-V_{L}}{R_{m,1}}\\
C_{2}\frac{dV_{2}}{dt} & =\frac{V_{1}-V_{2}}{R_{a,2}}+\frac{V_{3}-V_{2}}{R_{a,3}}-\frac{V_{2}-V_{L}}{R_{m,2}}\\
C_{3}\frac{dV_{3}}{dt} & =-g_{Ca}s^{2}r(V_{3}-V_{Ca})-g_{KCa}c_{Ca}(V_{3}-V_{K})+\frac{V_{2}-V_{3}}{R_{a,3}}-\frac{V_{3}-V_{L}}{R_{m,3}}\\
\frac{dc_{Ca}}{dt} & =g_{Ca}s^{2}rB(V_{Ca}-V_{3})-\frac{c_{Ca}}{\tau}\\
\frac{dP}{dt} & =\alpha_{P}(V)(1-P)+\beta_{P}(V)P,\quad P\in\{n,m,h,r,s\}.
\end{align*}
The opening and closing rates are given in the following table.

\begin{center}
\begin{tabular}{c|c}
Opening rate & Closing rate\tabularnewline
\hline 
$\alpha_{h}=70\exp(-50(V_{1}+0.07))$ & $\beta_{h}=\frac{1000}{1+\exp{(-100(V_{1}+0.0400)})}$\tabularnewline
$\alpha_{m}=10^{3}\psi(-100(V_{1}+0.045))$ & $\beta_{m}=4000\exp(-(V_{1}+0.07)/0.018)$\tabularnewline
$\alpha_{n}=100\psi(-100(V_{1}+0.06))$ & $\beta_{n}=125\exp(-12.5(V_{1}+0.07))$\tabularnewline
$\alpha_{r}=\begin{cases}
5, & \text{ if }V_{3}\leq-0.07\\
5\exp(-50(V_{3}+0.07)), & \text{ if }V_{3}>-0.07
\end{cases}$ & $\beta_{r}=5-\alpha_{r}$\tabularnewline
$\alpha_{s}=\frac{1600}{1+\exp(-72(V_{3}+0.005))}$ & $\beta_{s}=100\psi(200(V_{3}+0.0189))$\tabularnewline
\end{tabular}
\par\end{center}

The parameters are given by
\begin{gather*}
C_{1}=3.6\times10^{-11}F,\quad C_{2}=2\times10^{-11}F,\quad C_{3}=9.6\times10^{-15},\\
R_{m,1}=8.333\times10^{8}\Omega,\quad R_{m,2}=1.5\times10^{9}\Omega,\quad R_{m,3}=3.125\times10^{12}\Omega,\\
R_{a,2}=5\times10^{8}\Omega,\quad R_{a,3}=3\times10^{7}\Omega,\\
V_{Na}=0.045V,\quad V_{K}=-0.085V,\quad V_{Ca}=0.07V,\quad V_{L}=-0.0594V,\\
g_{Na}=5.4\times10^{-7}S,\quad g_{K}=5.4\times10^{-8}S,\quad g_{Ca}=9.6\times10^{-13}S,\quad g_{KCa}=7.68\times10^{-12},\\
I(t)=0.09\times10^{-9}A,\quad\tau=0.1s,\quad B=4.51389\times10^{12}.
\end{gather*}
The system has been solved on the interval $[0,0.1]$using the initial
values
\begin{gather*}
V_{1}(0)=0.07V,\quad V_{2}(0)=0.06V,\quad V_{3}(0)=0.06V,\\
c_{Ca}(0)=1.6\times10^{-4}\text{mol/m}^{3},\\
n(0)=0.8,\quad m(0)=1,\quad h(0)=0.3,\quad r(0)=1,\quad s(0)=0.11.
\end{gather*}

Figure~\ref{fig:Solution-of-maya} shows a plot of the solution.
\end{description}
All methods have been implemented in Matlab.\footnote{Matlab release 2016a, The MathWorks, Inc., Natick, MA, USA.}
The experiments have been carried out by running the codes with varying
tolerances $\textrm{TOL}=10^{-2-k/8}$ for $k=0,\ldots,48$. In the
applications, the coarser tolerances are of most interest. The tolerance
requirements in the codes implementing the new method are modeled
according to those used in the Matlab ode suite: For $z=(x,y)$, the
criterion
\[
|\textrm{err}_{i}|\leq\textrm{TOL}|z_{i}|+\textrm{AbsTol}_{i}
\]
 shall be satisfied for all components $z_{i}$ of $z$. Here, $\textrm{err}_{i}$
denotes the error estimate for $z_{i}$. $\textrm{AbsTol}_{i}$ has
been chosen as a typical size of $|z_{i}|$ multiplied by TOL.

The diagrams contain the obtained accuracy versus the computational
effort. The accuracy is measured as the error of the numerical solution
at the final time. Since the analytical solutions are not known, the
systems have been solved with very tight tolerances using \texttt{ode15s}
in order to obtain a reference solution.

The computational effort has been measured in terms of function evaluations.
One function evaluation corresponds to a computation of both $f$
and $g$. A Jacobian computations is counted as expensive as one function
evaluation. So one step of Hines' method corresponds to two function
evaluations while one step of the modified method needs 2.5 function
evaluations. This is consistent with the statistics provided by the
codes of the matlab ode suite.

\begin{figure}
\centering{}\includegraphics[height=0.8\textwidth]{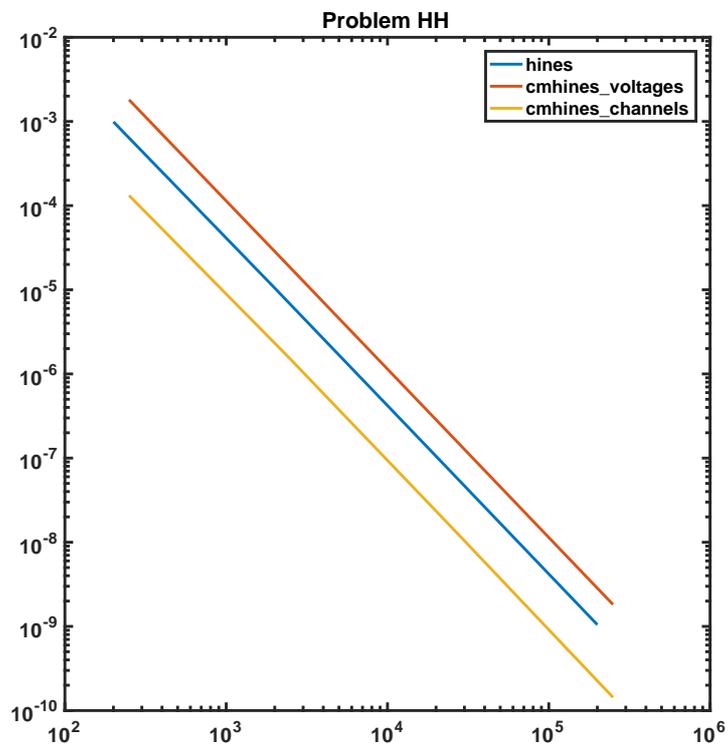}\caption{The Hodkin-Huxley system: Constant step size solvers\label{fig:HH-const}}
\end{figure}

\begin{figure}
\centering{}\includegraphics[height=0.8\textwidth]{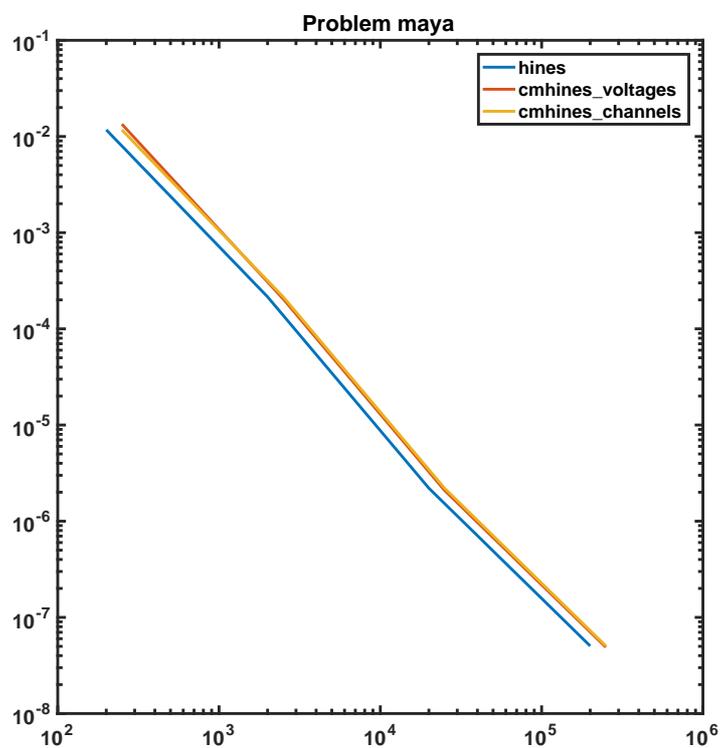}\caption{The soma-dendrite-spine system: Constant step size solvers\label{fig:maya-const}}
\end{figure}

\begin{sidewaysfigure}
\begin{centering}
\begin{tabular}{cc}
\includegraphics[width=0.4\textwidth]{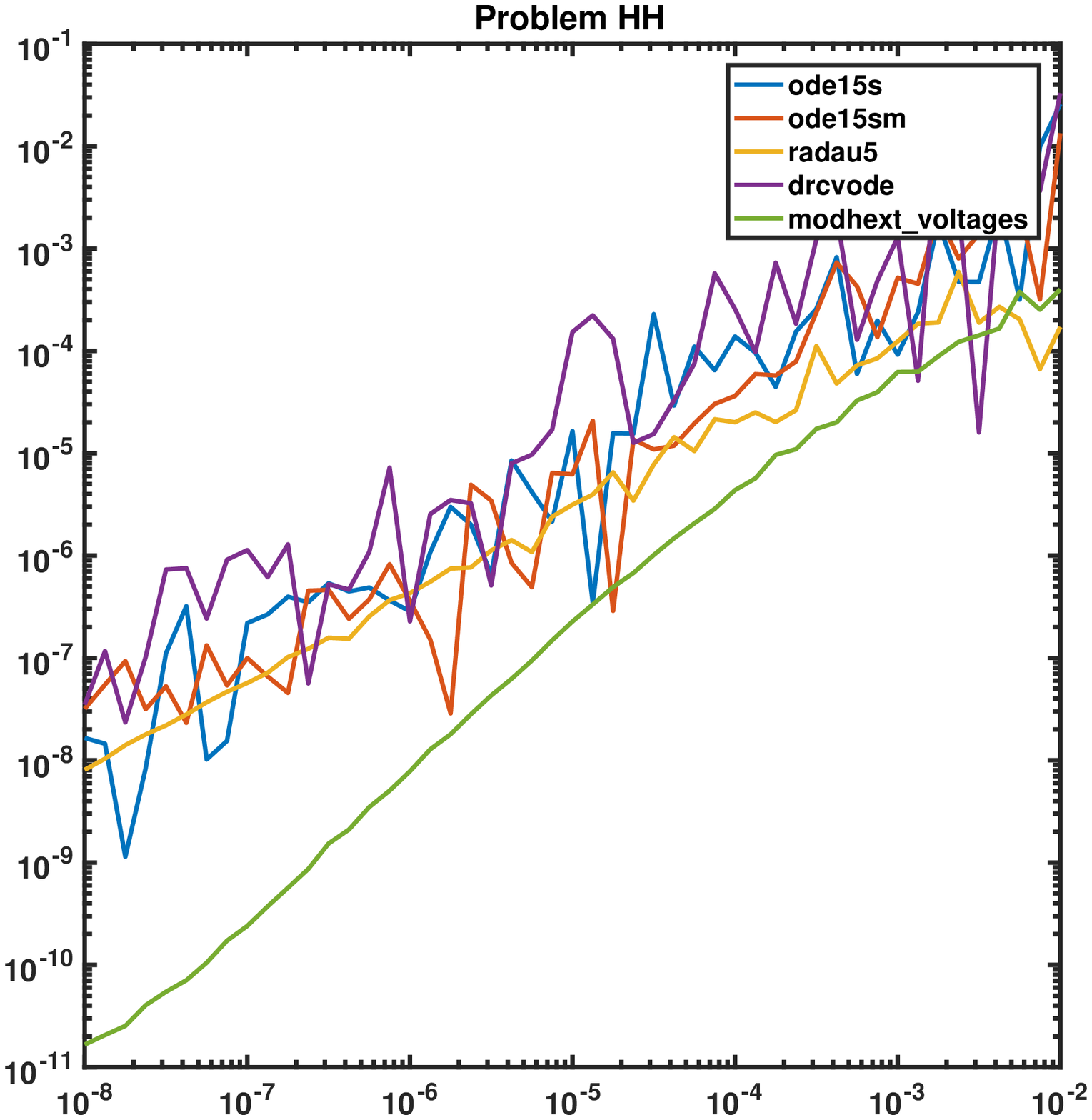} & \includegraphics[width=0.4\textwidth]{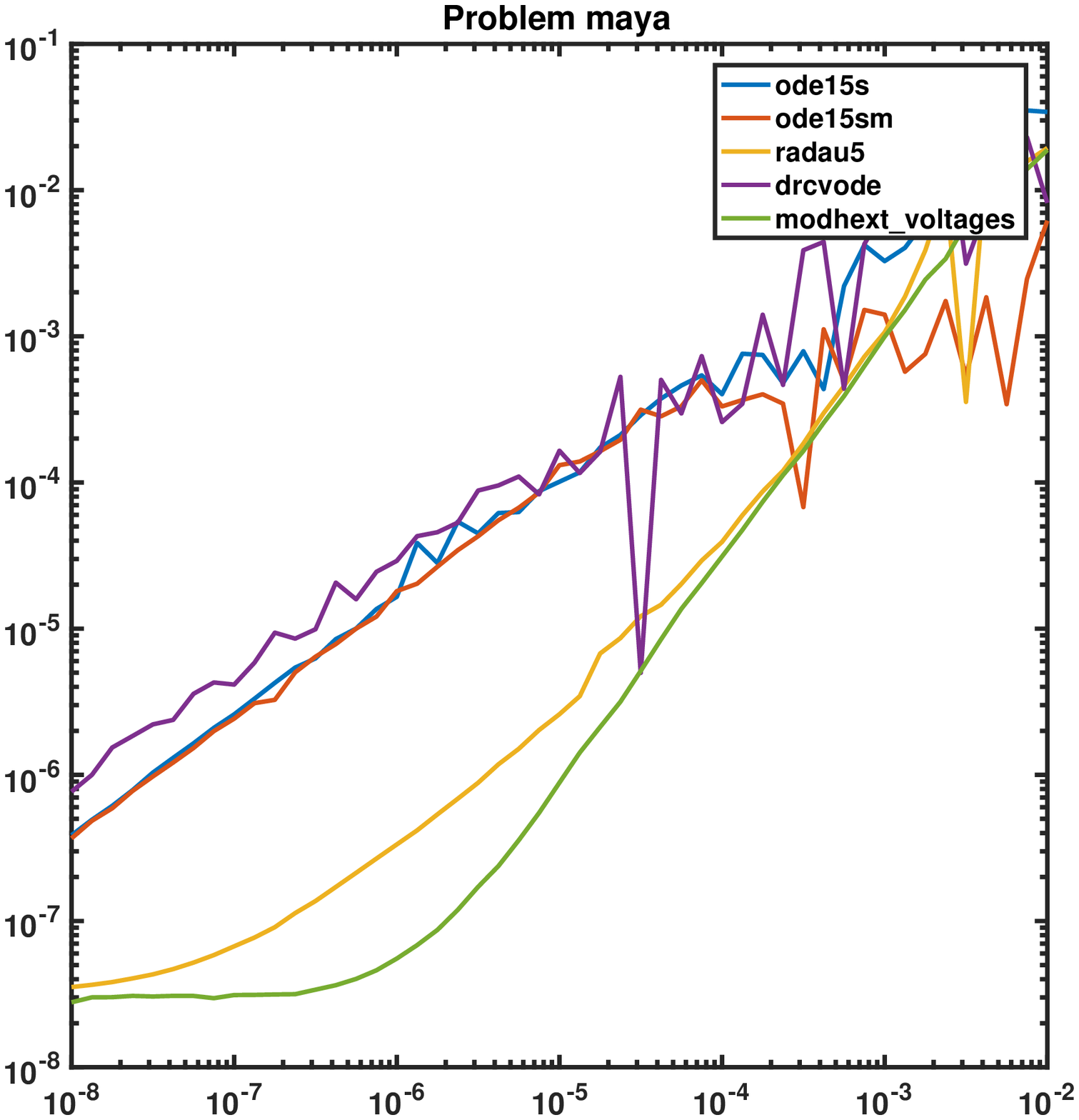}\tabularnewline
(a) & (b)\tabularnewline
\end{tabular}
\par\end{centering}
\caption{Tolerance-accuracy plots. (a) Hodgkin-Huxley model, (b) soma-dendrite-spine
system\label{fig:Tol-acc}}
\end{sidewaysfigure}

The following codes have been compared:
\begin{description}
\item [{\texttt{hines,~cmhines}}] This are the implementations of the
original method (\ref{eq:hinesa}) \textendash{} (\ref{eq:hinesb})
and the modification (\ref{eq:soeda}) \textendash{} (\ref{eq:soedc}),
respectively, for constant step sizes. The main purpose of using these
codes consists of showing second order convergence of both as well
as comparing the relative accuracy.
\item [{\texttt{modhines,~modhext,~modhnew}}] These are the new implementations
with error estimation and step size control according to the descriptions
above. Since these methods are no longer symmetric in $x$ and $y$
(in contrast to the original method), the experiments are run in two
versions: one where the voltages are used as $x$-component, and one
where the gate variables of the channels are used as $x$-components.
\item [{\texttt{ode15s,~ode23s,~ode23t,~ode23tb}}] This are Matlab's
ode solvers. Here, the system is solved as a whole. In these codes
the complete Jacobian is used. As an experiment, we modified even
\texttt{ode15s} in such a way that only the diagonal blocks of the
Jacobian are used (\texttt{ode15sm}). We intended to understand how
important the off-diagonal blocks are in the given examples. It should
be noted, however, that this approach is questionable since this may
break internal control strategies.
\item [{\texttt{radau5}}] This is the Fortran code RADAU5 developed by
Hairer\&Wanner \cite[Section IV.8]{HaWa96}. We used the interface
in \cite{Lu13} to call this code from Matlab.
\item [{\texttt{drcvode}}] This is the cvodes code from the Sundial package
\cite{HiBrGrLeSeShWo05}, version 2.8.2, using the matlab interface
version 2.5.0 \cite{Se05}.
\end{description}
The following observations can be made:
\begin{itemize}
\item Figures~\ref{fig:HH-const} and \ref{fig:maya-const} show clearly
that the new methods as well as Hines' method have second order of
accuracy. While Hines' method is symmetric in both the $x$- and $y$-components,
the symmetry is broken in the new method. Therefore, it becomes important
how the splitting is defined. For the Hodgkin-Huxley system, a much
better accuracy is obtained if the gates are chosen as $x$-components.
This difference is not seen in the soma-dendrite-spine example. Here,
the difference in efficiency between Hines' method and the new method
is mainly due to the fact that the new method is slightly more expensive
per step.
\item In both examples, the usage of higher order methods is preferred,
in particular in the case of higher accuracies. Thus, it is only the
version with local extrapolation which is competitive, at least for
low tolerances.
\item Among the second order methods, the additional flexibility offered
by variable step size solvers compared to the constant step size Hines'
method does not seem to pay off.
\item It is surprising how irregular the effort-accuracy curve is for many
of the well-established solvers. This holds in particular for the
Hodgkin-Huxley system. This is emphasized in Figure~\ref{fig:Tol-acc}
where the accuracy is plotted versus the tolerance requirement.
\end{itemize}

\section{Conclusions}

In the present note we have introduced a modification of a method
proposed by Hines for solving the large system of ordinary differential
equations arising when simulating large networks of neurons. The basic
motivation for the new method was to allow for a step size variation
while at the same time retaining the possibility for an efficient
linear algebra by using the special structure of the systems as it
is done in Hines' method. The relation of the new method to Strang
splitting and the Peaceman-Rachford method has been shown.

Since the systems are often stiff, an investigation of the domain
of absolute stability has been done. Here, a test equation inspired
by similar considerations for partitioned Runge-Kutta methods was
of much use. It turned out that, in many cases, these domains are
bounded.

The method has been implemented in different versions in matlab including
also an implementation of Richardson extrapolation. A number of tests
and comparisons to standard state-of-the-art solvers for ordinary
differential equations have been done. In the tests it turned out
that higher order methods are most efficient even in the case of rather
low tolerance which are of most practical interest.

The competivity of the new method compared to standard solvers depends
mostly on an efficient implementation of the linear algebra involved.
This could not be tested in the matlab environment. Therefore, we
will implement the new method in actual neuron simulators in the future
in order to obtain more realistic comparisons.

\bibliographystyle{plain}
\bibliography{P}

\clearpage{}

\pagebreak{}

\pagebreak{}

\pagebreak{}

\begin{sidewaysfigure}
\centering{}%
\begin{tabular}{cc}
\includegraphics[width=0.4\textwidth]{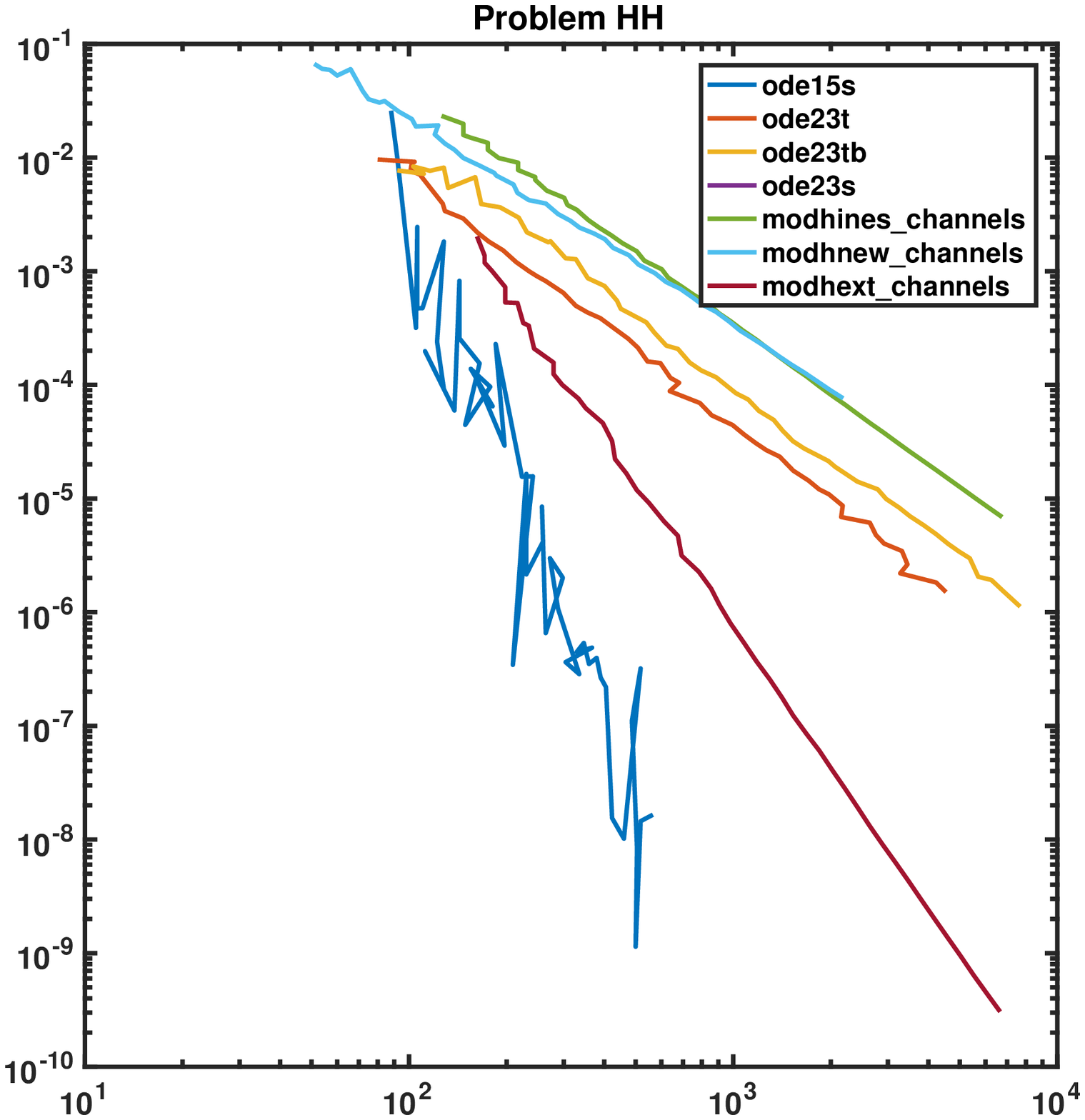} & \includegraphics[width=0.4\textwidth]{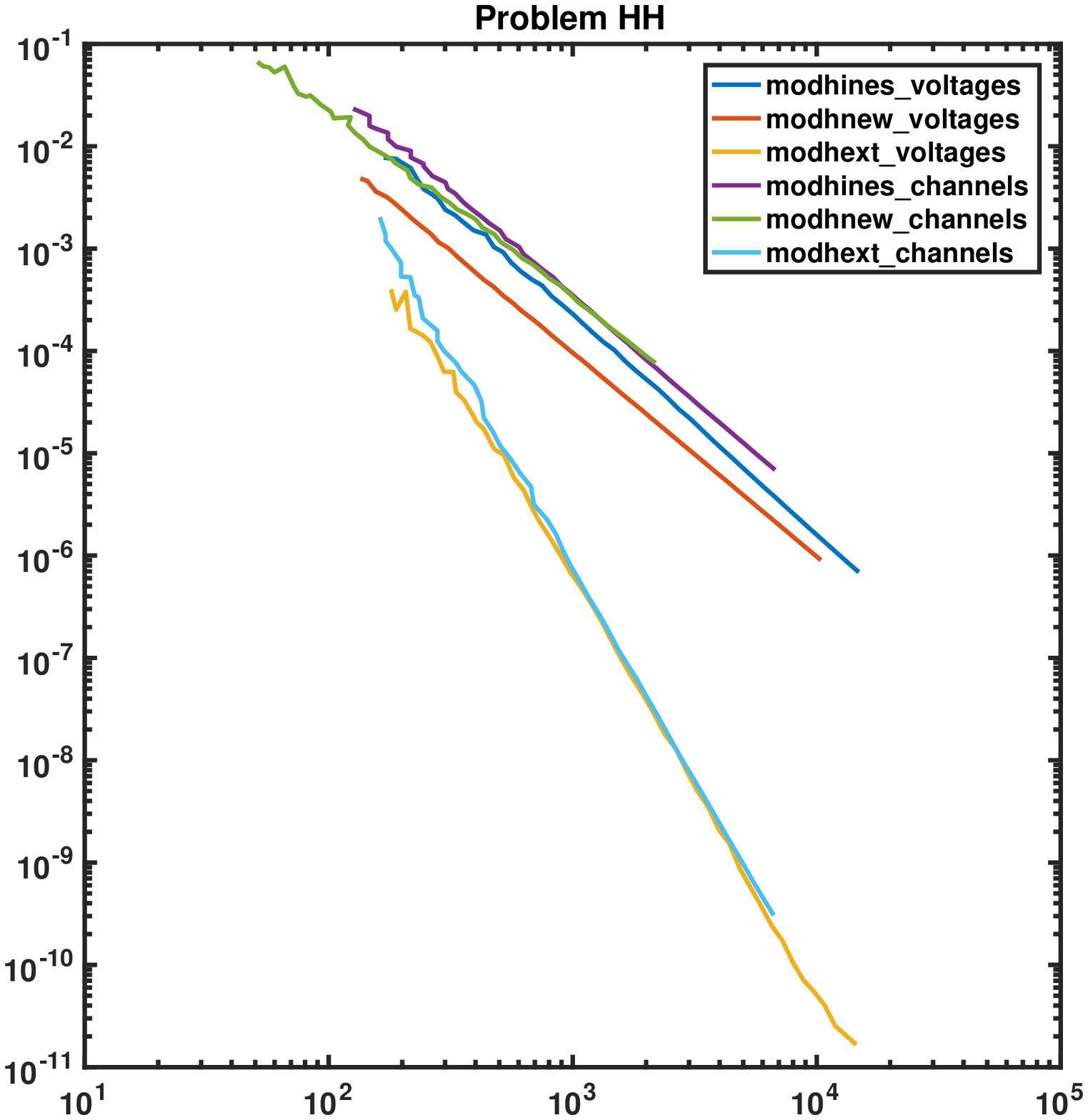}\tabularnewline
(a) & (b)\tabularnewline
\end{tabular}\caption{The Hodkin-Huxley system: Comparison of solvers. (a) behavior of the
new solvers in comparison to matlab's ode solvers, (b) behavior of
the new solvers with different assignments to the $x$- and $y$-components.
The tags ``voltages'' and ``channels'' indicate which set of variables
has been considered as $x$-components}
\end{sidewaysfigure}

\pagebreak{}
\begin{sidewaysfigure}
\centering{}%
\begin{tabular}{cc}
\includegraphics[width=0.4\textwidth]{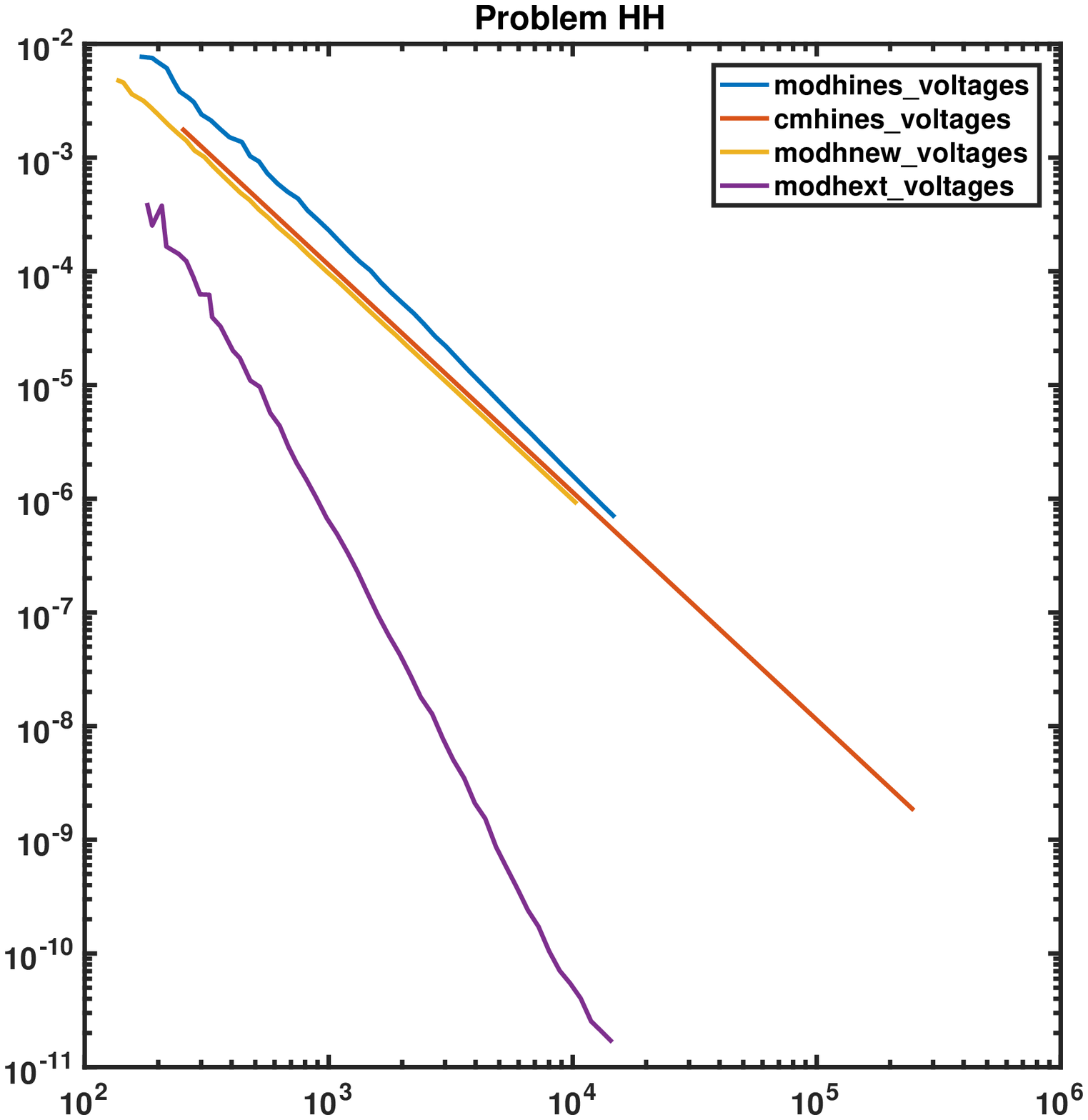} & \includegraphics[width=0.4\textwidth]{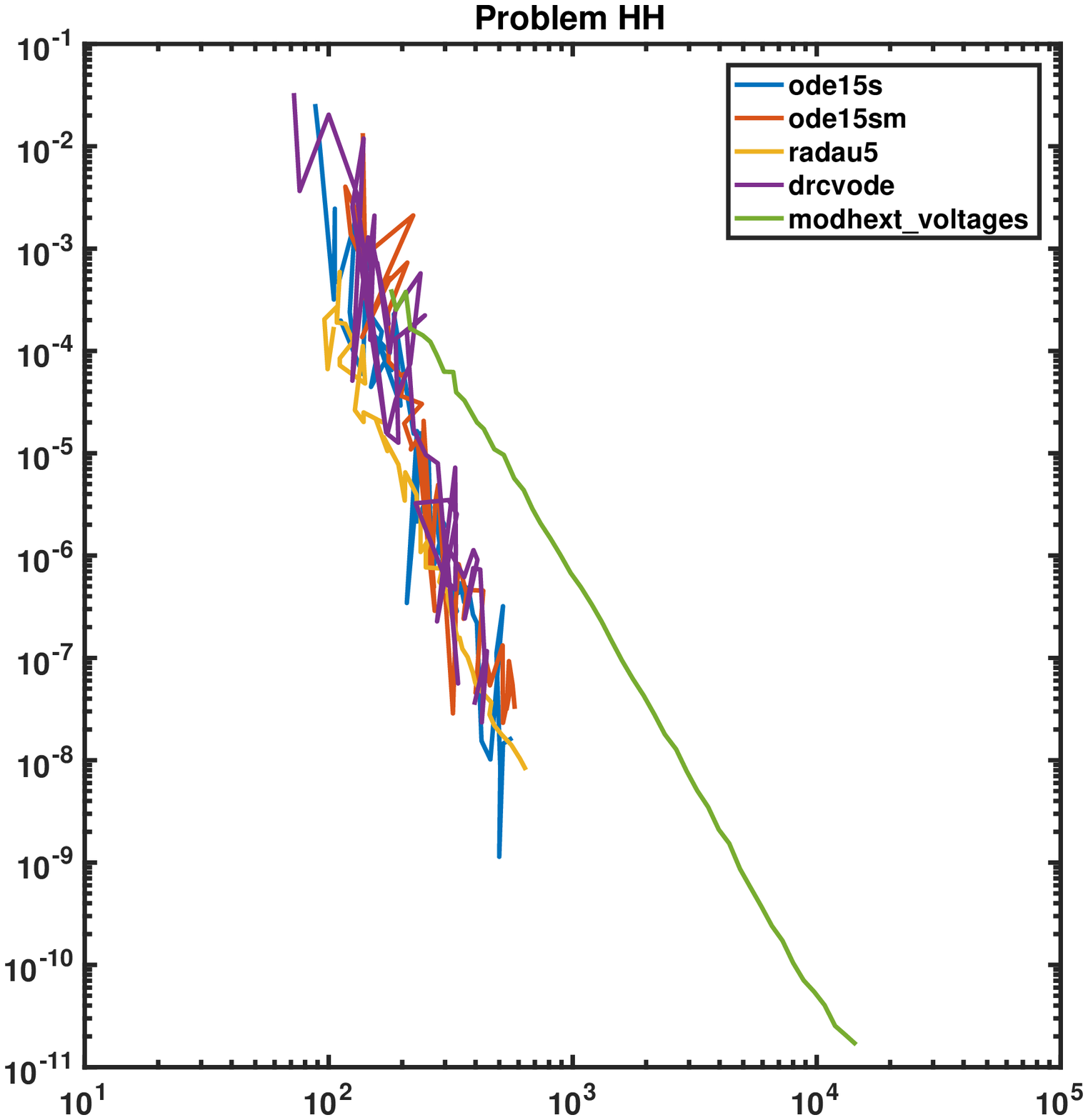}\tabularnewline
(a) & (b)\tabularnewline
\end{tabular}\caption{The Hodkin-Huxley system: Comparison of solvers. (a) The new solvers
and the constant step size Hines' method, (b) state-of-the-art methods
and the most efficient new version}
\end{sidewaysfigure}

\pagebreak{}

\pagebreak{}

\begin{sidewaysfigure}
\centering{}%
\begin{tabular}{cc}
\includegraphics[width=0.4\textwidth]{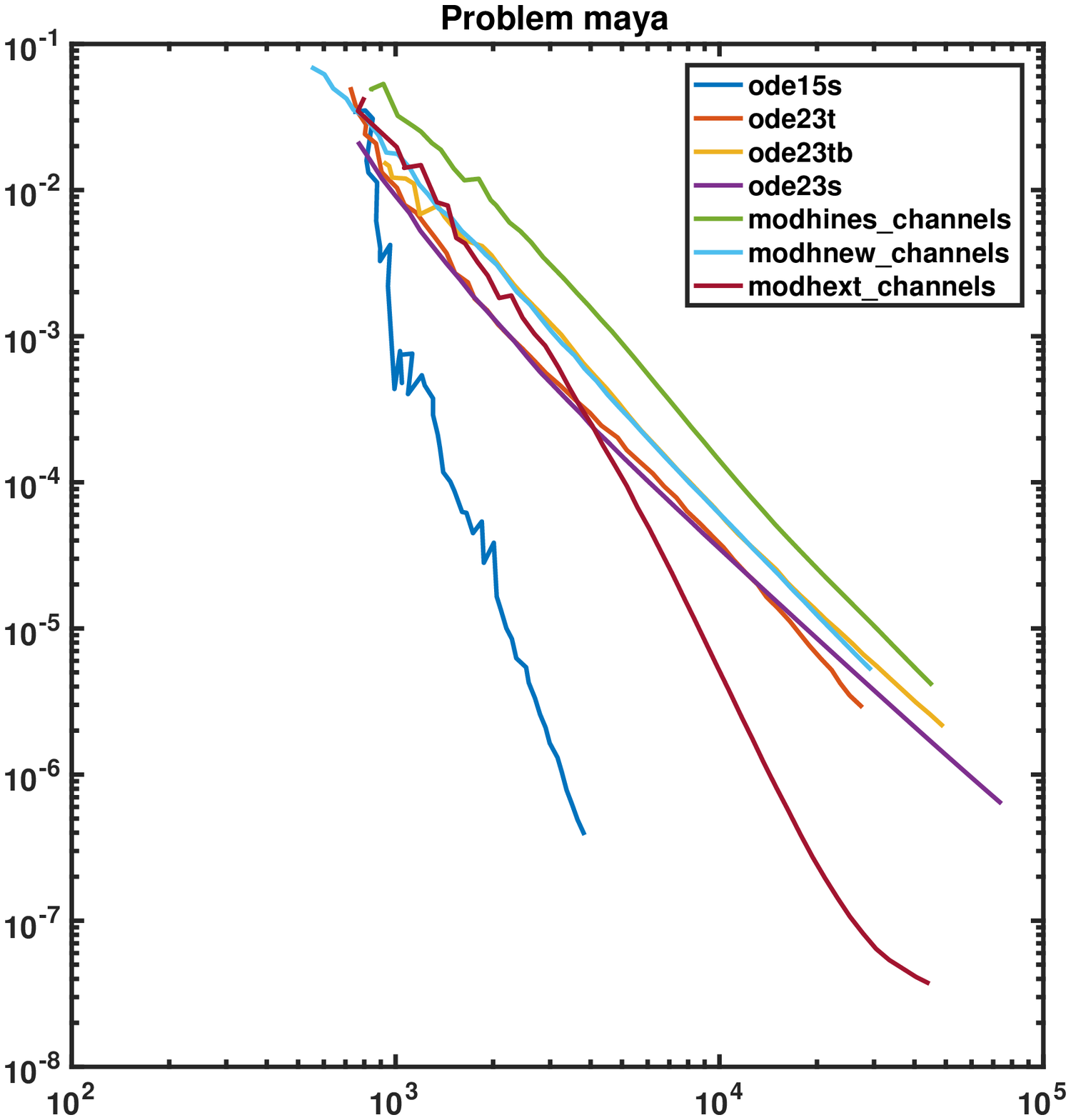} & \includegraphics[width=0.4\textwidth]{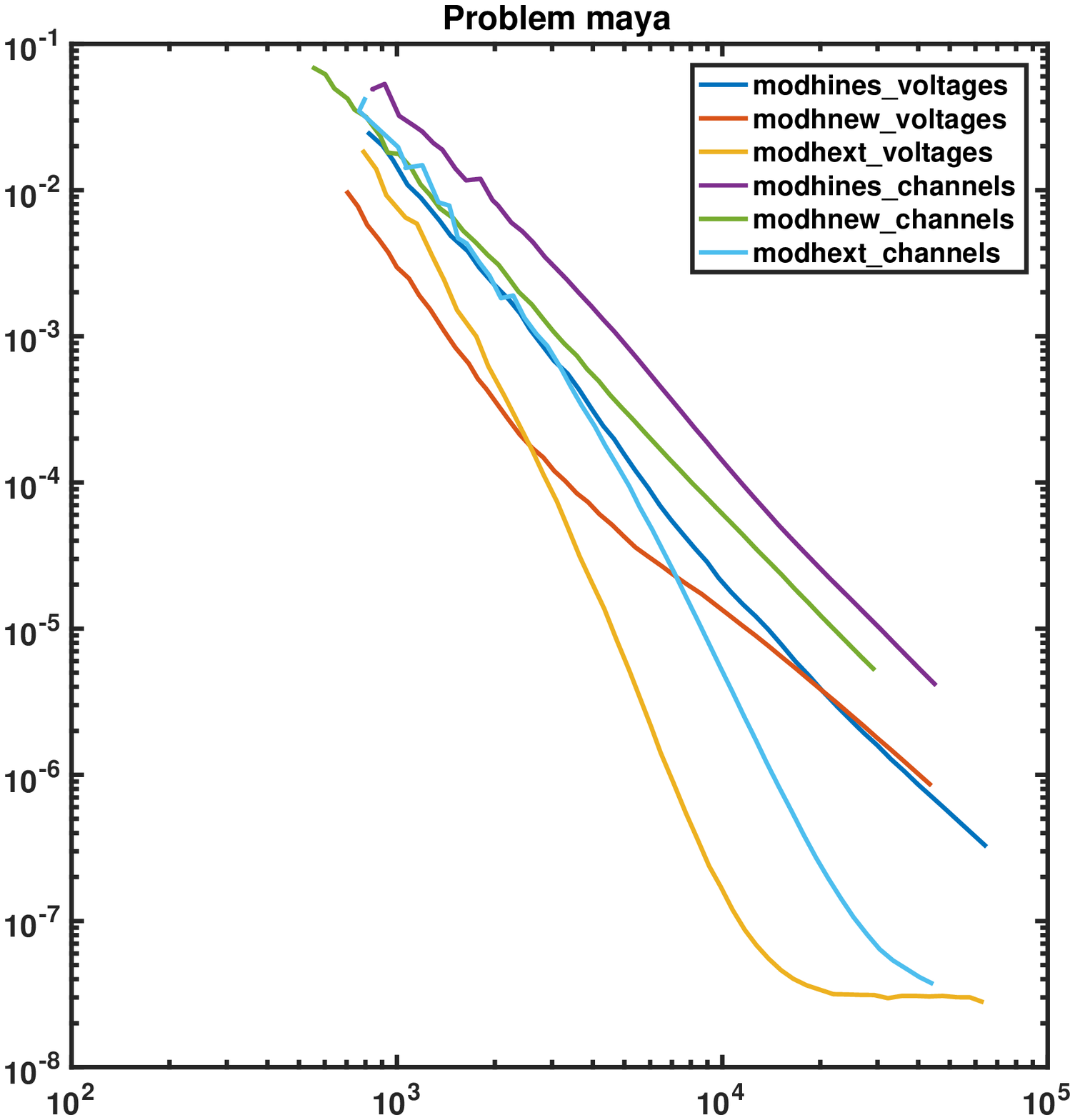}\tabularnewline
(a) & (b)\tabularnewline
\end{tabular}\caption{The soma-dendrite-spine system: Comparison of solvers. (a) behavior
of the new solvers in comparison to matlab's ode solvers, (b) behavior
of the new solvers with different assignments to the $x$- and $y$-components.
The tags ``voltages'' and ``channels'' indicate which set of variables
has been considered as $x$-components}
\end{sidewaysfigure}

\pagebreak{}
\begin{sidewaysfigure}
\centering{}%
\begin{tabular}{cc}
\includegraphics[width=0.4\textwidth]{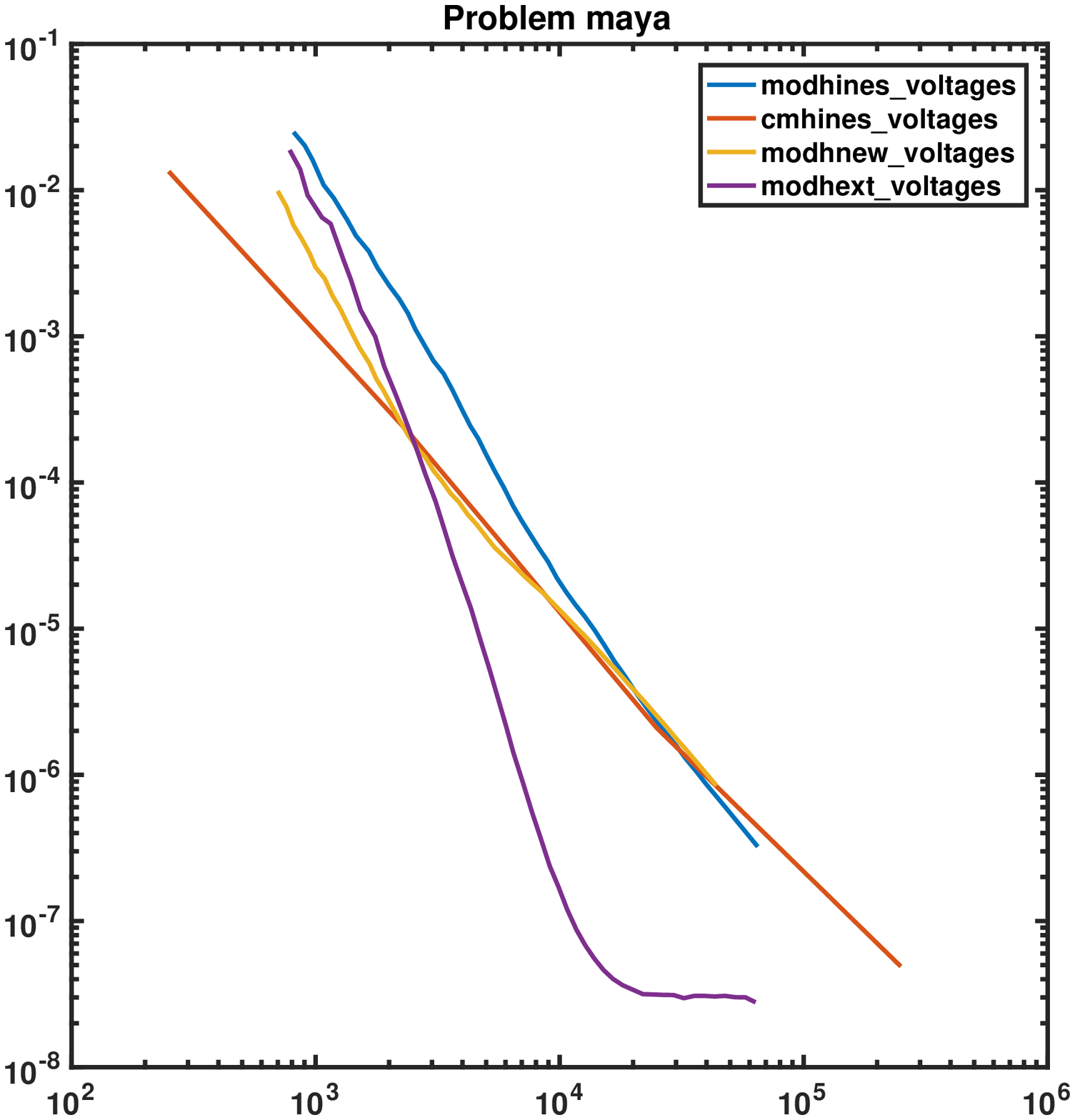} & \includegraphics[width=0.4\textwidth]{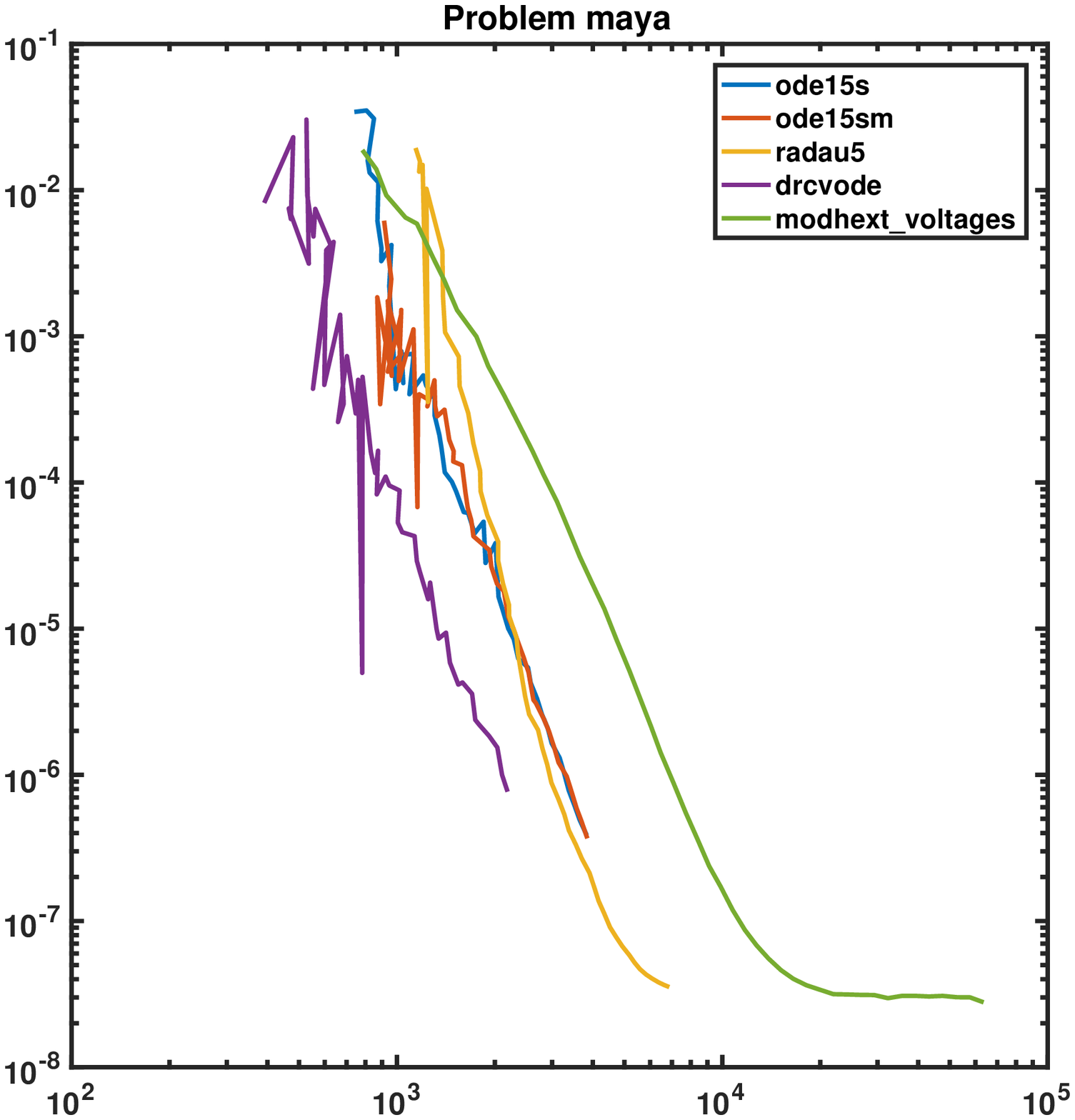}\tabularnewline
(a) & (b)\tabularnewline
\end{tabular}\caption{The soma-dendrite-spine system: Comparison of solvers. (a) The new
solvers and the constant step size Hines' method, (b) state-of-the-art
methods and the most efficient new version}
\end{sidewaysfigure}

\pagebreak{}

\end{document}